\documentclass[9pt]{amsart}
\usepackage{latexsym,amssymb,amsmath,mathabx}
\usepackage{graphicx}
\usepackage{color}
\usepackage[emtex,matrix,arrow]{xy}
\newtheorem{thm}{Theorem}
\newtheorem{cor}[thm]{Corollary}
\newtheorem{ex}[thm]{Example}

\newtheorem{lem}[thm]{Lemma}
\newtheorem{prop}[thm]{Proposition}
\newtheorem{defn}[thm]{Definition}
\newtheorem{rem}[thm]{Remark}
\numberwithin{equation}{section}

\newcommand{\Lb}{\mathbb L}

\newcommand{\Natural}{\mathbb N}

\newcommand{\To}{\rightarrow}

\newcommand{\Mm}{\mathcal{M}}
\newcommand{\Hh}{\mathcal{H}}

\newcommand{\Ff}{\mathcal{F}}
\newcommand{\Gg}{\mathcal{G}}
\newcommand{\Ii}{\mathcal{I}}

\newcommand{\Ll}{\mathcal{L}}

\newcommand{\Cc}{\mathcal{C}}

\newcommand{\Oo}{\mathcal{O}}

\newcommand{\Ee}{\mathcal{E}}

\newcommand{\F}{\mathbb{F}}
\newcommand{\Fb}{\mathbf{F}}
\newcommand{\Vv}{\mathcal{V}}
\newcommand{\Rb}{\mathbf{R}}

\newcommand{\Rr}{\mathcal{R}}

\newcommand{\G}{\mathbb{G}}

\newcommand{\Cgg}{\mathbf{C}}
\newcommand{\Fgg}{\mathbf{F}}

\newcommand{\Dgg}{\mathbf{D}}

\newcommand{\ngg}{\mathfrak{n}}
\newcommand{\id}{\rm{id}}

\newcommand{\homs}{\mathcal{H}om}

\newcommand{\sets}{\mathcal{S}et}

\newcommand{\modb}{\mathcal{M}od}

\newcommand{\catg}{\mathbf{Cat}}

\newcommand{\rep}{\mathcal{R}ep}
\newcommand{\equivb}{\mathbb{E}q}
\newcommand{\repg}{\mathbf{Rep}}
\newcommand{\ev}{\mathcal{V}ect}
\newcommand{\dev}{\mathbf{2Vect}}
\newcommand{\yon}{{\sf Yon}}

\newcommand{\adcat}{\mathbf{AdCat}}
\newcommand{\pseq}{\mathcal{P}sEq}

\begin{document}

\title{On the regular representation of an (essentially) finite 2-group}
\author{Josep Elgueta }
\address{Dept. Matem\`atica Aplicada II \\ Universitat
  Polit\`ecnica de Catalunya}
\email{josep.elgueta@upc.edu}

\begin{abstract}
The regular representation of an essentially finite 2-group $\G$
in the 2-category $\dev_k$ of (Kapranov and Voevodsky) 2-vector
spaces is defined and cohomology invariants
classifying it computed. It is next shown that all hom-categories
in $\repg_{\dev_k}(\G)$ are 2-vector spaces under quite standard
assumptions on the field $k$, and a formula giving the
corresponding ``intertwining numbers'' is obtained which proves
they are symmetric. Finally, it is shown that the forgetful
2-functor $\mbox{\boldmath$\omega$}:\repg_{\dev_k}(\G)\To\dev_k$
is representable with the regular representation as representing
object. As a consequence we obtain a $k$-linear equivalence
between the 2-vector space $\ev_k^{\Gg}$ of functors from the
underlying groupoid of $\G$ to $\ev_k$, on the one hand, and the
$k$-linear category $\Ee nd(\mbox{\boldmath$\omega$})$ of
pseudonatural endomorphisms of $\mbox{\boldmath$\omega$}$, on the
other hand. We conclude that $\Ee
nd(\mbox{\boldmath$\omega$})$ is a 2-vector space, and we (partially) describe 
a basis of it.
\end{abstract}

\maketitle

\section{Introduction}

Representation theory of 2-groups, i.e. of categories with a
structure analogous to that of a group, is a quite recent subject.
Although the special case of {\it discrete} 2-groups (2-groups
whose underlying category is discrete) was already considered in
the 1990's as (weak) actions of groups on categories (see
\cite{pD97}), the first works concerning general 2-groups appeared
as preprints in the current decade (\cite{CY03}, \cite{FB03},
\cite{BM04}, \cite{jE4}).

But what is a representation of a 2-group? By a
representation of a group it is meant its representation as a
group of automorphisms of an object in some category, mostly the
category $\sets_f$ of finite sets or the category $\ev_k$ of (finite
dimensional) vector spaces over a 
field $k$. Similarly, by a representation of a 2-group $\G$ one
means its representation as a 2-group of (weak) automorphisms of
an object in some 2-category $\Cgg$. For instance, in a
representation of $\G$ in the 2-category $\catg$ of (small)
categories, functors and natural transformations the objects of
$\G$ are thought of as self-equivalences of a certain category
$\Cc$ and the morphisms as natural isomorphisms between these
self-equivalences. This considerably generalizes, for instance,
the theory of representations of (finite) groups as permutations
of a (finite) set, recovered as the representations of the
associated (finite) discrete 2-group when $\Cc$ is a (finite)
discrete category.

Clearly, the first question one has to face
when studying representations of 2-groups is what 2-category we
should take as $\Cgg$. In \cite{jE4} I considered representations
of 2-groups in the 2-category $\dev_k$ of Kapranov and Voevodsky
2-vector spaces over a field $k$. This is a higher dimensional
version of the category $\ev_k$ where the role of the field $k$ is
played by the (semiring) category $\ev_k$.

The natural question arises whether this is a good choice. The answer
obviously depends on what one 
means by ``good''. A reasonable measure of the ``goodness'' of a
representation theory seems to be the amount of information on the
2-group we are able to recover from the corresponding (2-)category
of representations. In the case of groups, a representation theory
which has proved good, at least for some kinds of groups, is the
theory of complex finite dimensional linear representations. Under
appropriate assumptions on the group, it can indeed be completely
recovered from the corresponding category of such representations.
Results of this kind are generically known as {\it reconstruction
theorems}. The first such theorem, going back to the 1930's, is
Pontryagin's duality theorem on the canonical isomorphism between
any locally compact abelian topological group and its topological
bidual \cite{lsP87}. In this case, we are able to recover the
original group from just the group of isomorpism classes of
1-dimensional representations. Later on, Tannaka and Krein
concentrated on the problem of reconstructing any compact
topological group $G$, not necessarily abelian, from the whole
ring of isomorphism classes of finite dimensional linear
representations. Stated in a more modern categorical language
\cite{SR72}, they proved that the canonical map $\pi:G\To{\rm
End}(\omega)$ sending any $g\in G$ to the endomorphism of the
forgetful functor $\omega:\rep_{\ev_k}(G)\To\ev_k$ with components
$\pi(g)_{(V,\rho)}=\rho(g)$ defines an isomorphism of topological
groups between $G$ and the group ${\rm
Aut}_{\otimes}(\omega)\subset{\rm End}(\omega)$ of {\it monoidal}
automorphisms of $\omega$.

Although $\dev_k$ is introduced as a sort of higher dimensional
analog of $\ev_k$ it is pretty clear that for many 2-groups the
representation theory in $\dev_k$ will have deficiencies. Indeed,
it is easy to see \cite{jE4} that a representation of a 2-group
$\G$ in $\dev_k$ is given, among other things, by a representation
of the group $\pi_0(\G)$ of isomorphism classes of objects of $\G$
as automorphisms of a {\it finite} set. Here we think of
$\pi_0(\G)$ as a group with the group law induced by the product
existing between objects. Hence, for infinite 2-groups, in
particular, for `Lie 2-groups' (see \cite{BL03}) there will be
very few such representations and we will not be able to
reconstruct the whole 2-group from them.

However, it is plausible that this representation theory is good
enough if we restrict to essentially finite 2-groups, i.e.
2-groups whose underlying category has a finite set of isomorphism
classes of objects and a finite set of morphisms between any given
objects.

This paper is a natural continuation of the research program
initiated in \cite{jE4} with the purpose of investigating the representation
theory of 2-groups in this kind of 2-vector spaces. More particularly, it
arises as a first step 
toward the proof of the previous guess. Indeed one of the goals of the
program is to prove that any essentially finite 2-group $\G$ can
be recovered as the 2-group $\mathbb{A}{\rm
ut}_{\otimes}(\mbox{\boldmath$\omega$})$ of monoidal automorphisms
of the forgetful 2-functor
$\mbox{\boldmath$\omega$}:\repg_{\dev_k}(\G)\To\dev_k$ mapping any
representation of $\G$ in $\dev_k$ to its underlying 2-vector
space. This would translate into the category setting the
classical result we have for finite groups and its finite
dimensional linear representations.

Indeed, for a finite group $G$ the above mentioned theorem
identifying $G$ with the group ${\rm Aut}_{\otimes}(\omega)$ can
be proved using the regular representation $L(G)$ of $G$ and the
fundamental fact that this representation represents the forgetful
functor $\omega:\rep_{\ev_k}(G)\To\ev_k$. By the (enriched version
of the) Yoneda lemma it follows that there exists a linear
isomorphism $\phi:L(G)\stackrel{\cong}{\To}{\rm End}(\omega)$,
which is essentially an extension of the canonical map
$\pi:G\To{\rm End}(\omega)$ mentioned before. The point is that
$L(G)$ has a structure of a Hopf algebra whose group-like part is,
on the one hand, isomorphic to $G$ and, on the other hand,
bijectively mapped by $\phi$ to the subset of monoidal
endomorphisms of $\omega$.

With this situation in mind, the main purpose of this work is to
introduce an analogue of the regular representation for
essentially finite 2-groups $\G$ and to see, using the appropriate
2-categorical version of the Yoneda lemma, that it indeed
represents the corresponding forgetful 2-functor
$\mbox{\boldmath$\omega$}$. For this to make sense, it is first
necessary to prove that the 2-category of representations of an
essentially finite 2-group in $\dev_k$ is `closed' in the sense
that all its hom-categories are still 2-vector spaces. This is not
true for an arbitrary field $k$, but we shall prove it under quite
standard assumptions on $k$. This allows us to define a $k$-linear equivalence
of categories
$\ev_k^{\Gg}\stackrel{\simeq}{\longrightarrow}\mathcal{E}nd(\mbox{\boldmath$\omega$})$
analogous to the above $k$-linear isomorphism
$\phi:L(G)\stackrel{\cong}{\To}{\rm End}(\omega)$, where $\ev_k^{\Gg}$ denotes
the $k$-linear category of all $\ev_k$-valued functors on the underlying
groupoid $\Gg$ of $\G$. In a future paper it is intended to prove that both
categories actually admit a natural structure of a {\it Hopf 2-algebra} (higher
dimensional analog of a Hopf algebra) and that this equivelence is in fact as
Hopf 2-algebras, providing again an analog in our category setting of well known
results in the context of groups. 

\vspace{0.5 truecm}
\noindent{\bf Outline of the paper.} The first three sections
serve to recall some definitions and known facts needed in the
sequel. Specifically, Section 2 contains a quick review on
2-groups, including their description up to the relevant notion of
equivalence, and the basic definitions concerning the
representation theory of 2-groups. In section 3 we recall the
notion of Kapranov and Voevodsky 2-vector space, give some
examples (in particular, the 2-vector space underlying the regular
representation of an essentially finite 2-group) and discuss the
`closedness' of the corresponding 2-category. The classification
of the (general linear) 2-group of self-equivalences of an
arbitrary 2-vector space is also recalled here. Finally, in
Section 4 we recall from \cite{jE4} the cohomological description
of the representations of a 2-group in $\dev_k$.

The core of the
paper starts with Section 5, where we define the regular
representation of an essentially finite 2-group and explicitly
compute a set of data which classifies it up to equivalence
(Proposition~\ref{classificacio_repr_regular}).

In Section 6 it is
shown that, under appropriate assumptions, the 2-category of
representations of an essentially finite 2-group $\G$ in $\dev_k$
is indeed `closed' in the above sense. The main result is Theorem
21, where it is shown that all hom-categories are equivalent to a
product of categories of projective representations (with given
central charges) of a certain family of subgroups of $\pi_0(\G)$. We also
obtain a formula for 
computing the ranks of the 2-vector spaces one obtains as
categories of intertwiners, analogous to the so
called intertwining numbers, and we show that they are symmetric.

Finally, in
Section 7 we prove that 
the regular representation of an essentially finite 2-group
represents the forgetful 2-functor by identifying a universal
object in the underlying 2-vector space of the representation
(Theorem~\ref{functor_universal}). This allows us to obtain the above
mentioned $k$-linear equivalence between this 2-vector space and the
category $\Ee nd(\mbox{\boldmath$\omega$})$ of (weak)
endomorphisms of the forgetful 2-functor
$\mbox{\boldmath$\omega$}$, and to identify a `basis' of $\Ee
nd(\mbox{\boldmath$\omega$})$. Since any
$k$-linear functor on $\Ee nd(\mbox{\boldmath$\omega$})$ is determined, up to
isomorphism, but the image of a basis, having available a basis may be useful in
defining more structure on $\Ee nd(\mbox{\boldmath$\omega$})$, such as a
product or a coproduct. These are expected to play an important role in the
proof of the above mentioned reconstruction of $\G$ as a 2-group of symmetries
of $\mbox{\boldmath$\omega$}$.

\vspace{0.5 truecm}
\noindent {\bf Notation and terminology.} All over the
paper $k$ denotes a fixed field and $k^*=k\setminus\{0\}$. When we write
2-{\it something} 
we always mean the strict version. Sometimes, this is emphasized
by writting explicitly the word {\it strict}. The only exception
to this rule is when $something=group$, in which case we always
mean the weak version in general. Strict 2-groups are named so.
Vertical and horizontal compositions of natural transformations
and more generally, of 2-morphisms $\tau,\sigma$ in any 2-category
are respectively denoted by $\tau\cdot\sigma$ and
$\tau\circ\sigma$. For any set $X$ (resp. category $\Cc$), $X[0]$
(resp. $\Cc[0]$) denotes the corresponding discrete category with
only identity arrows (resp. locally discrete 2-category with only
identity 2-arrows). For any monoid $M$ (resp. monoidal category
$\Mm$), $M[1]$ (resp. $\Mm[1]$) denotes the corresponding
one-object category (resp. one-object 2-category). For any natural
number $n\geq 1$, $[n]$ denotes the set $\{1,\ldots,n\}$. $\ev_k$
denotes the category of finite dimensional vector spaces over $k$.

\vspace{0.5 truecm} {\bf Acknowledgements.} I would like to thank
Bertrand Toen for the many conversations we had during my stay at
the Laboratoire \'Emile Picard (Universit\'e Paul Sabatier), where
I started thinking about the subject of this paper.

\section{Review on 2-groups and their 2-categories of representations}

We assume the reader is familiar with the basic notions on
bicategories and in particular, with their one-object versions, the monoidal categories. See for instance \cite{jB67} or the short account
\cite{tL98}.

\subsection{Quick review on 2-groups}
\label{2grups} By a {\it 2-group} or {\it categorical group} it is
meant a monoidal groupoid $\G=(\Gg,\otimes,I,a,l,r)$ such that
each object $A$ has a weak inverse, i.e. an object $A^*$ such that
$A\otimes A^*\cong I\cong A^*\otimes A$. When the monoidal
groupoid is strict (the associator $a$ and the left and right unit
constraints $l,r$ are identities) and all inverses $A^*$ are
strict ($A\otimes A^*=I=A^*\otimes A$) the 2-group is called {\it
strict}.

The simplest examples are groups $G$ thought of as discrete
categories $G[0]$, and abelian groups $A$ thought of as one-object
2-groups $A[1]$. In both cases, the tensor porudct is given by the group
law. More generally, for any $G$-module $A$ we have 
the so called {\it split 2-group} $A[1]\rtimes G[0]$. Its set of objects is
$G$, its set of morphisms is $A\times G$, with a pair $(a,g)$
being an automorphism of $g$, and composition and tensor product
are given by
\begin{align*}
(a',g)\circ(a,g)&=(a'+a,g),
\\ g\otimes g'&=gg', \\ (a,g)\otimes(a',g')&=(a+g\cdot a',gg').
\end{align*}
This is just a special case of the general notion of semidirect product
for 2-groups, in this case between $G[0]$ and $A[1]$ (see \cite{GI01}).

In general, 2-groups arise as symmetries of objects in a
2-category. Thus for any 2-category $\Cgg$ and any object $X$ of
$\Cgg$ the groupoid $\mathcal{E}q(X)$ of self-equivalences of $X$
and 2-isomorphisms between these has a canonical structure of a
2-group with the tensor product given by composition of
self-equivalences and horizontal composition of 2-morphisms. We
shall denote by $\equivb(X)$ the 2-group so defined. Notice that
it is strict as a monoidal groupoid because $\Cgg$ is assumed to
be strict. However, $\equivb(X)$ is a non-strict 2-group in
general because there may exist objects having no strict inverse
(not all self-equivalences of $X$ will be isomorphisms).

As expected, 2-groups are the objects of a 2-category ${\bf 2Grp}$
whose 1-morphisms are monoidal functors between the corresponding
monoidal groupoids. Hence these are given by pairs $\F=(F,\mu)$
with $F:\Gg\To\Gg'$ a functor and $\mu$ a collection of natural
isomorphisms $\mu_{A,B}:F(A\otimes B)\stackrel{\cong}{\To}
F(A)\otimes' F(B)$ indexed by pairs of objects of $\Gg$ and
satisfying suitable coherence conditions. As it concerns
2-morphisms, they are given by the so called monoidal natural
transformations between these monoidal functors. See \cite{sM98}
for the precise definitions.

A basic result about 2-groups, due to Sinh \cite{hxSi75}, says
that any 2-group $\G$ is equivalent to a sort of ``twisted''
version of a split 2-group $A[1]\rtimes G[0]$ for some $G$-module
$A$. More precisely, let $\pi_0(\G)$ be the group of isomorphism
classes of objects in $\G$ with the product induced by the tensor
product, and let $\pi_1(\G)$ be the abelian group ${\rm Aut}(I)$
of automorphisms of the unit object of $\G$. This indeed is an
abelian group and it has a canonical $\pi_0(\G)$-module structure.
Then Sinh's classification theorem says that $\G$ is equivalent to
the semidirect product $\pi_1(\G)[1]\rtimes\pi_0(\G)[0]$ but
equipped with a non-trivial associator $a_{g,g',g''}:gg'g''\To
gg'g''$ given by
$$
a_{g,g',g''}=(\alpha(g,g',g''),gg'g''),
$$
where $\alpha\in Z^3(\pi_0(\G),\pi_1(\G))$ is a certain 3-cocycle
somehow constructed from the associator of $\G$. We shall denote
the 2-group defined in this way by
$\pi_1(\G)[1]\rtimes_{\alpha}\pi_0(\G)[0]$. For more details cf.
\cite{hxSi75} or the more accessible reference \cite{BL03}. The
groups $\pi_0(\G)$ and $\pi_1(\G)$ are called the {\it homotopy
groups} of $\G$ and the cohomology class $[\alpha]\in
H^3(\pi_0(\G),\pi_1(\G))$ its {\it Postnikov invariant}. Thus
split 2-groups are those whose Postnikov invariant is
$[\alpha]=0$. Any 3-cocycle $\alpha$ in the Postnikov invariant of
$\G$ is called a {\it classifying 3-cocycle} of $\G$.

In this paper we will mainly concentrate on {\it essentially
finite} 2-groups, by which we mean 2-groups $\G$ both of whose
homotopy groups $\pi_0(\G)$ and $\pi_1(\G)$ are finite.

\subsection{Representation bicategories of a 2-group}
\label{2-categoria_representacions}

The category of representations of a group $G$ in a category
$\Cc$, such as $\ev_k$, is nothing but the functor category $\Ff
un(G[1],\Cc)$. Indeed, a functor $F:G[1]\To\Cc$ is given by an
arbitrary object $X$ of $\Cc$ and a morphism of groups
$\rho:G\To{\rm Aut}_{\Cc}(X)$, and it is easy to check that
morphisms between representations correspond to natural
transformation between the respective functors.

By analogy, for any bicategory (resp. 2-category) $\Cgg$ and any
2-group $\G$ the bicategory (resp. 2-category) of representations
of $\G$ in $\Cgg$ is defined as the pseudofunctor bicategory
(resp. 2-category)
$$
\repg_{\Cgg}(\G):={\bf PsFun}(\G[1],\Cgg).
$$
Hence objects are pseudofunctors $\Fgg:\G[1]\To\Cgg$, 1-morphisms
are pseudonatural transformations between these and 2-morphisms
are modifications of pseudonatural transformations. When the
notions of pseudofunctor and pseudonatural transformation are
unpacked we get for the objects and morphisms in
$\repg_{\Cgg}(\G)$ the same kind of things that we get for the
objects and morphisms in $\rep_{\Cc}(G)$. Thus a representation of
$\G$ in $\Cgg$ is given by a pair $\Fgg=(X,\F)$, with $X$ an
object of $\Cgg$ and $\F:\G\To\equivb(X)$ a morphism of 2-groups,
and a 1-morphism or {\it intertwiner} $\xi:(X,\F)\To(X',\F')$ is
given by a pair $\xi=(f,\Phi)$, with $f:X\To X'$ a 1-morphism in
$\Cgg$ and $\Phi$
a family of 2-isomorphisms in $\Cgg$
\begin{equation} \label{condicio_intertwiner}
\xymatrix{ \ar @{} [drr] |{\mbox{{\Large
        $\stackrel{\Phi(A)}{\Leftrightarrow}$}}}
X\ar[d]_{f}\ar[rr]^{F(A)} & & X\ar[d]^{f} \\ X'\ar[rr]_{F'(A)} & &
X'
}
\end{equation}
indexed by the objects $A$ of $\Gg$. These 2-isomorphisms come
from the weakening of the action preserving condition in the usual
notion of intertwiner. They have to be natural in $A$ and to
satisfy some coherence conditions.

In our new setting, however, we further have morphisms between
intertwiners. More precisely, given intertwiners
$(f,\Phi),(\tilde{f},\tilde{\Phi}):(X,\F)\To(X',\F')$ a 2-morphism
or {\it 2-intertwiner} between them is just a 2-morphism
$\tau:f\Rightarrow\tilde{f}$ in $\Cgg$ satisfying a naturality
condition which involves the 2-cells $\Phi(A)$ and
$\tilde{\Phi}(A)$. See \cite{jE4} for more details.

As in any bicategory, we also have a composition law between
intertwiners and two composition laws between 2-intertwiners.
Composition between intertwiners is given by the so called
``vertical composition'' of pseudonatural transformations. More
explicitly, if $\xi=(f,\Phi):(X,\F)\To(X',\F')$ and
$\xi'=(f',\Phi'):(X',\F')\To(X'',\F'')$ the composite
$\xi'\circ\xi$ is described by the pair $(f'f,\Phi'\ast\Phi)$,
with the 2-cell $(\Phi'\ast\Phi)(A)$ given by the pasting
\begin{equation*}
\xymatrix{ \ar @{} [ddrr] |{\mbox{{\Large
        $\stackrel{\Phi(A)}{\Leftrightarrow}$}}}
X\ar[dd]_{f}\ar[rr]^{F(A)} & & X\ar[dd]^{f} \\ & & \\
\ar @{} [ddrr] |{\mbox{{\Large
        $\stackrel{\Phi'(A)}{\Leftrightarrow}$}}}
X'\ar[dd]_{f'}\ar[rr]^{F'(A)} & & X'\ar[dd]^{f'} \\ & & \\ X''\ar[rr]_{F''(A)} & & X''
}
\end{equation*}
Notice that such a pasting only makes sense when $\Cgg$ is a
(strict) 2-category, as it is the case in what follows. Otherwise,
we should also include the appropriate associativity constraint
2-cells. As for the two compositions between 2-intertwiners, they
are given by the vertical and horizontal composition of the
corresponding 2-morphisms in $\Cgg$ .

\section{2-vector spaces.}
\label{2espais_vectorials}

\subsection{Definition and examples}
\label{definicio_2ev}

There exists various notions of 2-vector space. See \cite{BC03},
\cite{jE5}, \cite{KV94}, \cite{mN97}. In this work we shall use
the notion originally introduced by Kapranov and Voevodsky in
\cite{KV94} although in a different guise.

According to Kapranov and Voevodsky, a 2-vector space is just a
special kind of what they call a $\ev_k$-module category. Roughly,
this is a symmetric monoidal category $\Vv$, analogous to the
abelian group in a vector space, together with a functor
$\odot:\ev_k\times\Vv\To\Vv$, called the {\it action} of $\ev_k$
on $\Vv$, and suitable natural isomorphisms coming from the
weakening of the usual axioms for a multiplication by scalars.
Then a 2-vector space is defined as a $\ev_k$-module category
equivalent to $\ev_k^n$ for some $n\geq 0$. Here $\ev_k^n$ is
assumed to be equipped with the $\ev_k$-action induced by the
usual tensor product of vector spaces, i.e.
$$
V\odot(V_1,\ldots,V_n)=(V\otimes V_1,\ldots,V\otimes V_n).
$$
Instead of this definition, however, we shall use the following
equivalent one. It provides an intrinsic characterization of
2-vector spaces and it is much easier to handle.
\begin{defn}\label{2-espai_vectorial}
A 2-vector space is a (small) $k$-additive category $\Vv$ which
admits a finite (possibly empty) basis of absolutely simple
objects.
\end{defn}
\noindent By a $k$-additive category it is meant a category
enriched over $\ev_k$ (not just over the category $\mathcal{A}b$
of abelian groups) and with zero object and all binary biproducts.
By an absolutely simple object in such a category it is meant an
object having no nonzero subobjects other than itself and such
that its vector space of endomorphisms is 1-dimensional. By a finite basis of
absolutely simple objects it is meant a 
finite set of absolutely simple objects $\{V_1,\ldots,V_n\}$ such
that any nonzero object is isomorphic to a {\it unique} finite
biproduct of them. Stated in this way, the definition is due to
Neuchl \cite{mN97}.

Notice that, in contrast to what happens in the case of vector
spaces, the basis of absolutely simple objects in a 2-vector space
is unique (up to isomorphism, of course). This has important consequences as it
concerns the representation theory of 2-groups on these 2-vector
spaces.

It readily follows from the above definition that the cartesian
product $\Vv\times\Vv'$ of two 2-vector spaces $\Vv,\Vv'$ is a new
2-vector space. A basis of absolutely simple objects is
$$\{(V_1,0'),\ldots,(V_n,0'),(0,V'_1),\ldots,(0,V'_{n'})\},
$$
where
$\{V_1,\ldots,V_n\}$ and $\{V'_1,\ldots,V'_{n'}\}$ are bases of
$\Vv$ and $\Vv'$, respectively.

\begin{ex}{\rm
The standard examples of 2-vector spaces are the product
categories $\ev_k^n$ for any $n\geq 0$. A basis of absolutely
simple objects is given by the objects
$\{(0,\ldots,\stackrel{i)}{k},\ldots,0),\ i\in[n]\}$. Any 2-vector
space $\Vv$ is actually equivalent to $\ev_k^n$ for some $n\geq
0$, called the {\it rank} of $\Vv$.}
\end{ex}

\begin{ex}\label{ex2}{\rm
Let $G$ be a finite group and $k$ an algebraically closed field
whose characteristic is zero or prime to the order of $G$. Then
the category $\rep_{\ev_k}(G)$ of finite dimensional $k$-linear
representations of $G$ is a 2-vector space of rank equal to the
number of conjugacy classes of $G$. A basis of absolutely simple
objects is given by any set of representatives of the equivalence
classes of irreducible representations. This example generalizes
to the case of projective representations with a given (arbitrary)
central charge and more generally, to finite dimensional modules
over an arbitrary semisimple $k$-algebra (see
\S~\ref{representacions_projectives} below).}
\end{ex}

\begin{ex} \label{vectG} {\rm
For any essentially finite 2-group $\G$ the category $\ev_k^{\
\Gg}$ of all functors $F:\Gg\To\ev_k$ and natural transformations
between them is a 2-vector space of rank
$$
{\rm rank}(\ev_k^{\ \Gg})=|\pi_0(\G)||\pi_1(\G)|.
$$
Indeed, for any 2-group, essentially finite or not, it always
happens that the automorphism group of any object $A$ of $\Gg$ is
isomorphic to $\pi_1(\G)$, even when the underlying groupoid $\Gg$
is non-connected. Thus we have an equivalence of categories
\begin{equation} \label{equivalencia}
\Gg\simeq\coprod_{g\in\pi_0(\G)}\pi_1(\G)[1],
\end{equation}
and therefore
$$
\ev_k^{\ \Gg}\simeq\ev_k\ ^{\displaystyle{\coprod_{g\in\pi_0(\G)}\pi_1(\G)[1]}}\cong\prod_{g\in\pi_0(\G)}\ev_k^{\pi_1(\G)[1]}=\prod_{g\in\pi_0(\G)}\rep_{\ev_k}(\pi_1(\G)).
$$
The claim follows now from the previous example and the fact that
$\pi_1(\G)$ is a finite abelian group. In particular, let
$\pi_1(\G)^*$ be the dual group of $\pi_1(\G)$, i.e. the
group of all group morphisms $\chi:\pi_1(\G)\To k^*$. Then a basis of
absolutely simple objects is given by the family of functors
$$\{\eta_{\chi,g}:\Gg\To\ev_k,\ \chi\in\pi_1(\G)^*,\
g\in\pi_0(G)\}
$$
defined on objects $A$ by
$$
\eta_{\chi,g}(A):=\left\{ \begin{array}{ll} k, & \mbox{if $A\in g$} \\ 0, & \mbox{otherwise},\end{array}\right.
$$
and on morphisms $\varphi:A\To B$, with $A,B\in g$, by
$$
\eta_{\chi,g}(\varphi)=\chi(h^{-1}_{A,B}(\varphi))\ {\id}_k.
$$
Here $h_{A,B}:\pi_1(\G)\To{\rm Hom}(A,B)$ for $A\cong B$ denote
isomorphisms we necessarily have to fix if we want to specify any
particular set of basic functors $\eta_{\chi,g}$. Thus although
for any object $A$ there is a canonical~\footnote{Actually we have
two such canonical isomorphisms, corresponding to the two
canonical morphisms ${\rm End}(I)\To{\rm End}(X)$ existing for any
object $X$ in any monoidal category $\Cc$. In case $\Cc$ is a
2-group these morphisms are isomorphisms; cf. \cite{SR72},
\S~1.3.3.3} isomorphism
$\gamma_A:\pi_1(\G)\stackrel{\cong}{\To}{\rm Aut}(A)$ there is no
canonical choice for the isomorphisms $\pi_1(\G)\cong{\rm
Hom}(A,B)$ when $A\cong B$ but $A\neq B$. Specifying such
isomorphisms is best done by choosing representatives
$A_1,\ldots,A_k$ in each isomorphism class $g\in\pi_0(\G)$, with
$A_1$ equal to the unit object $I$ of $\G$, together with
isomorphisms $\iota_A:A\To A_i$ between each object $A$ and its
representative $A_i$. Making these choices actually amounts to
fixing an equivalence of categories as in (\ref{equivalencia}).
Then an isomorphism $h_{A,B}$ is given by
$h^{-1}_{A,B}(\varphi)=\gamma^{-1}_{A_i}(\iota_{B}\varphi\iota_A^{-1})$.
Different choices lead to different isomorphisms $h_{A,B}$ and
hence, to different (but isomorphic) basic functors
$\eta_{\chi,g}$. To get the decomposition of an arbitrary functor
$\eta:\Gg\To\ev_k$ as a biproduct of the $\eta_{\chi,g}$ we just
need to take the restriction of $\eta$ to the various subgroupoids
${\rm Aut}(A_i)[1]$ and decompose them as a direct sum of irreps.
}
\end{ex}

Let $\catg_k$ be the 2-category of all (small) $k$-linear
categories, $k$-linear functors and natural transformations. Then
we denote by $\dev_k$ its full sub-2-category with objects all
2-vector spaces. Observe that we still have a third 2-category in
between them. Namely, the full sub-2-category $\adcat_k$ of
$\catg_k$ with objects all $k$-additive categories.

For any two objects $\Vv,\Vv'$ in $\dev_k$ the corresponding
hom-category is denoted by $\homs_k(\Vv,\Vv')$ instead of
$\homs_{\dev_k}(\Vv,\Vv')$. Observe that $\dev_k$ is a {\it
repletive} sub-2-category of $\catg_k$ in the sense that any
object of $\catg_k$ equivalent (in $\catg_k$) to a 2-vector space
is itself a 2-vector space. In fact, any $k$-linear equivalence
between 2-vector spaces maps a basis of absolutely simple objects
to a basis of the same kind.

\subsection{Hom-categories in $\dev_k$}
\label{hom-categories_a_dev} As in the vector spaces setting, all
hom-categories in $\dev_k$ are themselves 2-vector spaces. Because
of its importance we include here the proof of this elementary but
fundamental fact.

\begin{prop} \label{homs_interns}
Let $\Vv,\Vv'$ be any 2-vector spaces of ranks $n,n'$
respectively. Then $\homs_{k}(\Vv,\Vv')$ is a 2-vector space of
rank $nn'$.
\end{prop}
\begin{proof}
The category $\homs_{k}(\Vv,\Vv')$ has an obvious $k$-additive
structure, with the `zero functor' mapping all objects of
$\Vv$ to any fixed zero object of $\Vv'$ as a zero object
of $\homs_k(\Vv,\Vv')$, and with the biproduct $H\oplus\tilde{H}$ of
any pair $H,\tilde{H}:\Vv\To\Vv'$ of $k$-linear functors 
computed pointwise.

The existence of a finite basis follows from the general fact that, up to
isomorphism, a $k$-linear functor $H:\Vv\To\Vv'$ is completely
given by the corresponding {\it matrix of ranks} $R=(r_{i'i})\in{\rm
  Mat}_{n'\times n}(\Natural)$. By definition, it
is the matrix whose entries are uniquely determined by the
condition
$$
H(V_i)\cong\bigoplus_{i'=1}^{n'}r_{i'i}V'_{i'},\qquad i\in[n],
$$
where $\{V_1,\ldots,V_n\}$ and $\{V'_1,\ldots,V'_{n'}\}$ are bases
of absolutely simple objects of $\Vv$ and $\Vv'$, respectively.
The matrix of ranks of the biproduct of two functors corresponds to taking the
sum of the respective matrices of ranks. Hence a basis of
$\homs_k(\Vv,\Vv')$ is given by any representative set of
$k$-linear functors
$$
\{H_{i'i},\ (i',i)\in[n']\times[n]\}
$$
whose isomorphism classes are described by the unit matrices
(matrices having a unique nonzero entry equal to 1).

Once we have fixed biproduct functors in $\Vv$ and $\Vv'$, it is easy to see
that any morphism
$\tau:H\Rightarrow\tilde{H}$ in $\homs_k(\Vv,\Vv')$ is completely
given by its `basic components', i.e. the components
$$
\tau_{V_i}:\bigoplus_{i'=1}^{n'}r_{i'i}V'_{i'}\To\bigoplus_{i'=1}^{n'}\tilde{r}_{i'i}V'_{i'},\quad i=1,\ldots,n
$$
for a basis $\{V_1,\ldots,V_n\}$ of $\Vv$. Moreover, each of these
components $\tau_{V_i}$ is in turn described by a collection of
$n'$ arbitrary matrices $M_{i'i}\in{\rm
Mat}_{\tilde{r}_{i'i}\times r_{i'i}}(k)$, $i'=1,\ldots,n'$, giving
the morphism between the homologous ``isotypic'' pieces
$$
M_{i'i}:V_{i'}'\oplus\stackrel{r_{i'i}}{\cdots}\oplus V'_{i'}\To
V'_{i'}\oplus\stackrel{\tilde{r}_{i'i}}{\cdots}\oplus
V'_{i'},\quad i'=1,\ldots,n'
$$
(if both $r_{i'i},\tilde{r}_{i'i}\neq 0$; otherwise, they are
empty matrices). See \cite{jE3} for more details. In particular, any natural
endomorphism of a basic functor $H_{i',i}$ is completely given by an
(arbitrary) scalar $\lambda\in k$, and this shows they are indeed absolutely
simple. 
\end{proof}

\subsection{General linear 2-groups}
For any 2-vector space $\Vv$ we shall denote by $\G\Lb(\Vv)$ the
corresponding 2-group of ($k$-linear) self-equivalences, and by
$\Gg\Ll(\Vv)$ the underlying groupoid. These 2-groups $\G\Lb(\Vv)$
should be thought of as analogs in our category setting of the
usual general linear groups, and they will be called {\it general
linear 2-groups}. The underlying monoidal groupoids are always
strict because $\dev_k$ is a strict 2-category. However, they are
non-strict 2-groups in general because there may exist no strict
inverses for objects. If $n$ is the rank of $\Vv$, it may be shown
that $\G\Lb(\Vv)$ is a split 2-group with homotopy groups
\begin{align*}
\pi_0(\G\Lb(\Vv))&\cong S_n, \\ \pi_1(\G\Lb(\Vv))&\cong(k^*)^n
\end{align*}
and with the usual action of $S_n$ on $(k^*)^n$. For the details,
see for ex. \cite{jE5}, where these 2-groups are computed for a
more general kind of 2-vector spaces including those of Kapranov
and Voevodsky.

\section{Linear representations of a 2-group.}

\subsection{Description up to equivalence}
\label{classes_iso_repr}

Let $\repg_{\dev_k}(\G)$ be the 2-category of representations of
$\G$ in $\dev_k$. Thus an object is a pair $\Fgg=(\Vv,\F)$ with
$\Vv$ a 2-vector space and $\F=(F,\mu):\G\To\G\Lb(\Vv)$ a morphism
of 2-groups. The rank of $\Vv$ is called the {\it dimension} of
the representation.

As in any 2-category, two objects $\Fgg$ and
$\Fgg'$ are said to be {\it equivalent} when there exists an
equivalence between them, i.e. a {\it weakly} invertible
intertwiner between them. In \cite{jE4} it is shown that the
equivalence class of a representation is completely specified by a quadruple
$(n,\rho,\beta,c)$ with 
\begin{itemize}
\item
$n\geq 0$ a natural number,
\item
$\rho:\pi_0(\G)\To S_n$ a morphism of groups, where $S_n$ denotes the symmetric group on $n$ elements,
\item
$\beta:\pi_1(\G)\To(k^*)^n_{\rho}$ a morphism of
$\pi_0(\G)$-modules such that $[\beta_*(\alpha)]=0$ (in the
cohomology group $H^3(\pi_0(\G),(k^*)^n_{\rho})$), and
\item
$c\in C^2(\pi_0(\G),(k^*)^n_{\rho})$ a normalized 2-cochain such that $\partial c=\beta_*(\alpha)$.
\end{itemize}
Here $\alpha$ is any classifying 3-cocycle of $\G$, and
$(k^*)^n_{\rho}$ denotes the abelian group of $n$-tuples of
nonzero elements of $k$ with the $\pi_0(\G)$-module structure
induced by $\rho$ and the usual action of $S_n$ on $(k^*)^n$
$$
g(\lambda_1,\ldots,\lambda_n)=(\lambda_{\rho(g^{-1})(1)},\ldots,\lambda_{\rho(g^{-1})(1)}),\quad
g\in\pi_0(\G).
$$
Notice that this description is neither canonical nor faithful. It
is non-canonical because it depends on the specific 3-cocycle
$\alpha$ we choose to describe $\G$ up to equivalence. In
particular, the 2-cochain $c$ changes with $\alpha$. But it is
also non-faithful because different quadruples, even for a fixed
$\alpha$, can describe the same equivalence class of
representations. More precisely, the two quadruples
$(n,\rho,\beta,c),(n',\rho',\beta',c')$ specify the same
equivalence class of representations if and only if $n=n'$ and
there exists $\sigma\in S_n$ such that
$\rho'=\sigma\rho\sigma^{-1}$, $\beta'=\sigma\beta$ and
$[c']=[\sigma c]$.

A specific representation $\Fgg=(\Vv,\F)$ whose equivalence class
is described by the quadruple $(n,\rho,\beta,c)$ is the following:
\begin{itemize}
\item $\Vv=\ev_k^n$.
\item $\F=(F,\mu):\G\To\G\Lb(\ev_k^n)$ is the monoidal functor defined as follows:
\begin{itemize}
\item
it maps $A\in{\rm Obj}\Gg$ to the permutation functor
$$
F(A)\equiv P_{\rho[A]}:\ev_k^n\To\ev_k^n
$$
acting on objects $(V_1,\ldots,V_n)$ by
$$
P_{\rho[A]}(V_1,\ldots,V_n):=(V_{\rho[A](1)},\ldots,V_{\rho[A](n)})
$$
($[A]$ denotes the isomorphism class of $A$ and $\rho[A]$ its
image by $\rho$);
\item
it maps a morphism $\varphi:A\To B$ of $\Gg$ to the natural
automorphism
$$
F(\varphi):P_{\rho[A]}\Rightarrow P_{\rho[B]}
$$
(notice that $[B]=[A]$) whose {\it basic
components}~\footnote{Once we fix specific direct sum functors in
the codomain category $\ev_k^m$, any natural transformation
$\tau:H\Rightarrow H'$ between $k$-linear functors
$H,H':\ev_k^n\To\ev_k^m$ is completely determined by the ``basic''
components $\tau_{(0,\ldots,\stackrel{i)}{k},\ldots,0)}$ for all
$i=1,\ldots,n$. This fact was already mentioned before in the
proof of Proposition~\ref{homs_interns}. See for ex. \cite{jE3}
for more details.} are
$$
F(\varphi)_{(0,\ldots,\stackrel{i)}{k},\ldots,0)}:=(0,\ldots,\beta_{\rho[A](i)}(h_{A,B}^{-1}(\varphi))\ {\id}_K,\ldots,0)
$$
(the isomorphisms $h_{A,B}$ are defined in Example~\ref{vectG}
above);
\item
for any objects $A,B$ of $\Gg$ the natural isomorphism
$$
\mu_{A,B}:P_{\rho[A\otimes B]}\Rightarrow P_{\rho[A]}\circ P_{\rho[B]}
$$
(actually, an automorphism) giving the monoidal structure is that
whose basic components are
$$
(\mu_{A,B})_{(0,\ldots,\stackrel{i)}{k},\ldots,0)}:=(0,\ldots,c_{\rho[A\otimes B](i)}([A],[B])\ {\id}_K,\ldots,0)
$$
\end{itemize}
\end{itemize}
We shall denote the representation so defined~\footnote{Relative
to the direct sum functors fixed in each 2-vector space $\ev_k^n$,
$n\geq 1$.} by $\Fgg(n,\rho,\beta,c)$. In particular, we see that
$n$ gives the dimension of the representation, $\rho$ and $\beta$
give the action of the corresponding functor $F:\Gg\To\Gg\Ll(\Vv)$
on objects and morphisms, respectively, and $c$ gives the monoidal
structure.

The morphism $\beta$ admits the following alternative description.
The left action of $\pi_0(\G)$ on $\pi_1(\G)$ induces a left
action on $\pi_1(\G)^*$ given by
$$
(g\chi)(u)=\chi(g^{-1}u),\quad g\in\pi_0(\G),\ \chi\in\pi_1(\G)^*,\
u\in\pi_1(\G).
$$
For any natural number $n\geq 1$ and any morphism of groups
$\rho:\pi_0(\G)\To S_n$, let $[n]_{\rho}$ be the set
$[n]\equiv\{1,\ldots,n\}$ equipped with the $\pi_0(\G)$-set
structure induced by $\rho$. Then we have the following.

\begin{lem} \label{lema}
For any pair $(n,\rho)$ as above a morphism of $\pi_0(\G)$-modules
$\beta:\pi_1(\G)\To(k^*)^n_{\rho}$ is the same thing as a
$\pi_0(\G)$-equivariant map $\gamma:[n]_{\rho}\To\pi_1(\G)^*$.
\end{lem}
\begin{proof}
From any $\beta$ as in the statement we define a map $\gamma$ also
as in the statement by $\gamma(i)=\beta_i$, $i=1,\ldots,n$. It is
easy to check that this sets a bijection between both types of
maps.
\end{proof}
\noindent This is the same kind of things that Crane and Yetter
\cite{CY03} and Baez et al. \cite{BBFW09} obtain for the
representations of 2-groups in Yetter's measurable categories.

\subsection{Some examples of linear representations}
\begin{ex}{\rm
The 1-{\it dimensional trivial representation}, denoted by $\Ii$,
is defined by the pair $(\Vv,\F)$ with $\Vv=\ev_k$ and $\F$ equal
to the trivial strict morphism of 2-groups. It corresponds to
$n=1$ and $\beta$ and $c$ the respective constant maps equal to
1.}
\end{ex}

\begin{ex}{\rm
Any $z\in Z^2(\pi_0(\G),k^*)$ defines a 1-dimensional
representation of $\G$ where $\rho$ and $\beta$ are trivial, and
cohomologous cocycles define different but equivalent
representations. In fact, for discrete 2-groups $G[0]$ we get in
this way a canonical bijection between $H^2(G,k^*)$ and the set of
equivalence classes of its 1-dimensional representations in
$\dev_k$.}
\end{ex}

\begin{ex}{\rm
More generally, for any $n\geq 1$ and any $[z]\in
H^2(\pi_0(\G),(k^*)^n)$ we have an $n$-dimensional representation
whose corresponding functor $F:\Gg\To\Gg\Ll(\ev_k^{\Gg})$ is the
trivial one mapping any object to the identity of $\ev_k^{\Gg}$
but equipped with a non-trivial monoidal structure. These are
called {\it cocyclic representations} of $\G$. }
\end{ex}
\begin{ex}{\rm
Any permutation representation $\rho:\pi_0(\G)\To S_n$ induces an
$n$-dimensional representation of $\G$ whose corresponding functor
$F:\Gg\To\Gg\Ll(\ev_k^{\Gg})$ just permutes the objects according
to $\rho$. These are called {\it permutation representations} of
$\G$. Equivalent permutation representations $\rho$ of $\pi_0\G)$
give rise to equivalent permutation representations of $\G$. In
this way, the theory of permutation representations of $\pi_0(\G)$
embeds into the theory of representations of $\G$ in $\dev_k$ (for
a more precise statement, see Theorem~5.13 in \cite{jE4}). }
\end{ex}

Clearly, a generic linear representation of $\G$ is a sort of mixture of a cocyclic and a permutation representation.

\section{The regular representation of an essentially finite 2-group.}

Recall that the {\it regular representation} of a group $G$ is the
permutation representation of $G$ induced by the left action of
$G$ on itself by left translations. Equivalently, it is the
representation defined by the vector space $L(G)$ of all functions
$f:G\To k$ with (left) $G$-action given by $(gf)(h)=f(hg)$. In this section we
describe an analog of this 
representation for essentially finite 2-groups and a
quadruple $(n,\rho,\beta,c)$ which classifies it up to
equivalence.

\subsection{Definition of the regular representation}
Let $\G=(\Gg,\otimes,I,a,l,r)$ be an essentially finite 2-group. A
canonical representation $\Rb=(\Vv_{\Rb},\F_{\Rb})$ of $\G$ can be
obtained as follows. Take as $\Vv_{\Rb}$ the 2-vector space
$\ev_k^{\Gg}$ (cf. Example~\ref{vectG}), and as
$F_{\Rb}:\Gg\To\Gg\Ll(\ev_k^{\Gg})$ the functor which maps
$A\in{\rm Obj}\Gg$ to the $k$-linear self-equivalence
$F_{\Rb}(A):\ev_k^{\Gg}\To\ev_k^{\Gg}$ acting on objects
$\eta:\Gg\To\ev_k$ and morphisms $\tau:\eta\Rightarrow\eta'$ by
$$
F_{\Rb}(A)(\eta):=\eta\circ(-\otimes A),\quad F_{\Rb}(A)(\tau):=\tau\circ 1_{-\otimes A}.
$$
If $\varphi:A\To B$ is any morphism of $\Gg$, $F_{\Rb}(\varphi)$ is the natural transformation
$$
F_{\Rb}(\varphi):F_{\Rb}(A)\Rightarrow F_{\Rb}(B):\ev_k^{\Gg}\To\ev_k^{\Gg}
$$
whose $\eta$-component $F_{\Rb}(\varphi)_{\eta}:\eta\circ(-\otimes A)\Rightarrow\eta\circ(-\otimes B)$ is defined by
$$
F_{\Rb}(\varphi)_{\eta,C}:=\eta({\id}_C\otimes\varphi),\qquad C\in{\rm Obj}\Gg.
$$
The point is that the functor $F_{\Rb}$ so defined has a canonical
monoidal structure induced by the associativity constraints in
$\G$. More precisely, we have the following:

\begin{lem}
For any $B,C\in{\rm Obj}\Gg$ let $\mu_{B,C}:F_{\Rb}(B\otimes
C)\Rightarrow F_{\Rb}(B)\circ
F_{\Rb}(C):\ev_k^{\Gg}\To\ev_k^{\Gg}$ be the natural
transformation with components
$\mu_{B,C;\eta}:\eta\circ(-\otimes(B\otimes
C))\Rightarrow\eta\circ(-\otimes C)\circ(-\otimes B)$ given by
$$
\mu_{B,C;\eta}:=1_{\eta}\circ a_{-,B,C},\qquad \eta\in{\rm Obj}\ev_k^{\Gg},
$$
where $a_{-,B,C}:-\otimes(B\otimes C)\Rightarrow(-\otimes
C)\circ(-\otimes B)$ is the natural isomorphism defined by the
associativity constraints $a_{A,B,C}:A\otimes(B\otimes
C)\cong(A\otimes B)\otimes C$ of $\G$. Then $\mu_{B,C}$ is natural
in $B,C$ and the collection $\mu=\{\mu_{B,C}\}_{B,C}$ provides
$F_{\Rb}$ with a monoidal structure.
\end{lem}
\begin{proof}
Note first that the diagram
$$
\xymatrix{
\eta\circ(-\otimes(B\otimes C))\ar@{=>}[rr]^{1_{\eta}\circ a_{-,B,C}}\ar@{=>}[d]_{\tau\circ 1_{-\otimes(B\otimes C)}} & & \eta\circ(-\otimes C)\circ(-\otimes B)\ar@{=>}[d]^{\tau\circ 1_{(-\otimes C)(-\otimes B)}} \\ \eta'\circ(-\otimes(B\otimes C))\ar@{=>}[rr]_{1_{\eta'}\circ a_{-,B,C}} & & \eta'\circ(-\otimes C)\circ(-\otimes B)
}
$$
commutes for any $\tau:\eta\Rightarrow\eta'$ by the interchange
law, so that $\mu_{B,C;\eta}$ is indeed natural in $\eta$.
Naturality of $\mu_{B,C}$ in $B,C$ means the commutativity of the
diagram
$$
\xymatrix{
F_{\Rb}(B\otimes C)\ar@{=>}[rr]^{\mu_{B,C}}\ar@{=>}[d]_{F_{\Rb}(\varphi\otimes\psi)} & & F_{\Rb}(B)\circ F_{\Rb}(C)\ar@{=>}[d]^{F_{\Rb}(\varphi)\circ F_{\Rb}(\psi)} \\ F_{\Rb}(B'\otimes C')\ar@{=>}[rr]_{\mu_{B',C'}} & & F_{\Rb}(B')\circ F_{\Rb}(C')
}
$$
for all morphisms $\varphi:B\To B'$, $\psi:C\To C'$ in $\Gg$.
Taking components this amounts to the commutativity of the
diagrams
$$
\xymatrix{
\eta(A\otimes(B\otimes C))\ar[rr]^{\eta(a_{A,B,C})}\ar[d]_{\eta(id_A\otimes(\varphi\otimes\psi))} & & \eta((A\otimes B)\otimes C)\ar[d]^{\eta((id_A\otimes\varphi)\otimes\psi)} \\ \eta(A\otimes(B'\otimes C'))\ar[rr]_{\eta(a_{A,B',C'})} & & \eta((A\otimes B')\otimes C')
}
$$
for all $\eta:\Gg\To\ev_k$ and all $A\in{\rm Obj}\Gg$, and these
diagrams commute because $a_{A,B,C}$ is natural in $B,C$. Finally,
since the underlying monoidal groupoid of $\G\Lb(\ev_k^{\Gg})$ is
strict, the coherence condition on $\mu$ reduces to the
commutativity of the diagram
$$
\xymatrix{
F_{\Rb}(B\otimes(C\otimes D))\ar@2[rr]^{F_{\Rb}(a_{B,C,D})}\ar@2[dd]_{\mu_{B,C\otimes D}} & & F_{\Rb}((B\otimes C)\otimes D)\ar@2[d]^{\mu_{B\otimes C,D}} \\ & & F_{\Rb}(B\otimes C)\circ F_{\Rb}(D)\ar@2[d]^{\mu_{B,C}\circ 1_{F_{\Rb}(D)}} \\ F_{\Rb}(B)\circ F_{\Rb}(C\otimes D)\ar@2[rr]_{1_{F_{\Rb}(B)}\circ\mu_{C,D}} & & F_{\Rb}(B)\circ F_{\Rb}(C)\circ F_{\Rb}(D)
}
$$
for any objects $B,C,D$ of $\Gg$. Taking again components this
amounts to the commutativity of the diagrams
$$
\xymatrix{
\eta(A\otimes(B\otimes(C\otimes D)))\ar[rr]^{\eta(id_A\otimes a_{B,C,D})}\ar[dd]_{\eta(a_{A,B,C\otimes D})} & & \eta(A\otimes((B\otimes C)\otimes D))\ar[d]^{\eta(a_{A,B\otimes C,D})} \\ & & \eta((A\otimes (B\otimes C))\otimes D)\ar[d]^{\eta(a_{A,B,C}\otimes id_{D})} \\ \eta((A\otimes B)\otimes (C\otimes D))\ar[rr]_{\eta(a_{A\otimes B,C,D})} & & \eta(((A\otimes B)\otimes C)\otimes D)
}
$$
for any $\eta:\Gg\To\ev_k$ and any objects $A,B,C,D$ of $\Gg$, and
these diagrams commute by the pentagon axiom on the associativity
isomorphisms.
\end{proof}

\begin{defn}
For any essentially finite 2-group $\G$ the regular representation
of\ \ $\G$ is the representation $\Rb$ defined by the pair
$(\ev_k^{\Gg},\F_{\Rb})$ with $\F_{\Rb}=(F_{\Rb},\mu)$ the above
morphism of 2-groups.
\end{defn}

\begin{ex}{\rm 
For any finite group $G=\{g_1,\ldots,g_n\}$, the regular representation of
$G[0]$ is the strict monoidal functor $F_{\Rb}:G[0]\To\Gg\Ll(\ev_k^G)$
mapping $g\in G$ to the permutation functor $\ev_k^G\To\ev_k^G$ given by
$(V_{g_1},\ldots,V_{g_n})\mapsto(V_{g_1g},\ldots,V_{g_ng})$.
}
\end{ex}

\begin{ex}{\rm 
For any finite abelian group $A$, the regular representation of $A[1]$ is
(equivalent to) the strict monoidal functor $R_{\Rb}:A[1]\To\Gg\Ll(\Rr
ep_{\ev_k}(A))$ mapping the unique object to the identity functor and $a\in A$
to the natural automorphism $F_{\Rb}(a):id\Rightarrow id$ defined by
$F_{\Rb}(a)_{(V,\rho)}=\rho(a)$ for any representation $(V,\rho)$ of $A$
(observe that $\rho(a)$ indeed is an intertwiner from the representation
$(V,\rho)$ to itself because $A$ is abelian). Thus it essentially reduces to the
canonical morphism from $A$ into the center $Z(\Rr ep_{\ev_k}(A))$ 
of its category of linear representations.
}
\end{ex}

\subsection{Classification}
\label{classificacio} Let $p=|\pi_0(\G)|$, $q=|\pi_1(\G)|$. We
know from Example~\ref{vectG} that $\Rb$ has dimension
$n_{\Rb}=pq$. In this subsection we describe a particular triple
$(\rho_{\Rb},\beta_{\Rb},c_{\Rb})$ of the kind described in
\S~\ref{classes_iso_repr} that classifies $\Rb$. Recall that such
a triple is unique only ``up to conjugation''. In particular, it
depends on the choice of a representative  of the Postnikov
invariant of $\G$. Let us fix once and for all such a
representative $\alpha\in Z^3(\pi_0(\G),\pi_1(\G))$, that we can
assume normalized without loss of generality.

Before describing the triple $(\rho_{\Rb},\beta_{\Rb},c_{\Rb})$
let us first introduce some notation. Let us denote by
$S_p\times\stackrel{q)}{\cdots}\times S_p\hookrightarrow S_{pq}$
the embedding mapping the $i^{th}$-factor $S_p$ ($i=1,\ldots,q$)
to the subgroup of $S_{pq}$ leaving all $j\in[pq]$ invariant
except the elements $\{(i-1)p+1,\ldots,ip\}$, which are permuted
accordingly. In terms of permutation matrices, this means mapping
the permutation matrices $(P_1,\ldots,P_q)$ to the block diagonal
permutation matrix $P={\rm diag}(P_1,\ldots,P_q)$. For any
linearly ordered finite group $G=\{g_1<\ldots<g_r\}$ let us
further denote by $\kappa:G\To S_r$ the composite
$$
G\hookrightarrow{\rm Aut}(G)\stackrel{\cong}{\To}S_r,
$$
where $G\hookrightarrow{\rm Aut}(G)$ denotes Cayley's embedding
mapping $g\in G$ to the right translation $g'\mapsto g'g^{-1}$,
and $\cong$ stands for the isomorphism of groups induced by the
chosen linear order in $G$.

The starting point to classify {\bf R} is the classification of
the general linear 2-groups $\G\Lb(\Vv)$ described in
\S~\ref{2espais_vectorials}. We know that $\pi_0(\ev_k^{\Gg})\cong
S_{pq}$, but we need to specify a particular such isomorphism. To
do this we choose a linear order in one of the sets of basic
functors $\{\eta_{\chi,g}\}$ for $\ev_k^{\Gg}$ described in
Example~\ref{vectG}. As explained before, we have various such
sets of basic functors and we fix any one of them. Let us further
fix linear orders $g_1<\ldots<g_p$ in $\pi_0(\G)$ and
$\chi_1<\ldots<\chi_q$ in $\pi_1(\G)$, and take as linear order in
the fixed set of basic functors the lexicographical one, i.e.
$\eta_{\chi_1,g_1}<\ldots<\eta_{\chi_1,g_p}<\ldots<\eta_{\chi_q,g_1}<\ldots<\eta_{\chi_q,g_p}$.
This way a permutation $\sigma\in S_{pq}$ becomes identified with
the isomorphism class of the corresponding permutation functor
$\ev_k^{\Gg}\To\ev_k^{\Gg}$. Moreover, this automatically
specifies a particular isomorphism
$\pi_1(\ev_k^{\Gg})\cong(k^*)^{pq}$, namely that sending
$u:{\id}_{\ev_k^{\Gg}}\Rightarrow{\id}_{\ev_k^{\Gg}}$ to the
corresponding basic components
$(u_{\chi_1,g_1},\ldots,u_{\chi_1,g_p},\ldots,u_{\chi_q,g_1},\ldots,u_{\chi_q,g_p})$,
which we know are completely given by one non-zero scalar each of
them (see proof of Proposition~\ref{homs_interns}). With these
choices we have the following.

\begin{prop} \label{classificacio_repr_regular}
The equivalence class of $\Rb$ is described by the following triple $(\rho_{\Rb},\beta_{\Rb},c_{\Rb})$:
\begin{itemize}
\item[(i)]
$\rho_{\Rb}:\pi_0(\G)\To S_{pq}$ is given by the composite
\begin{equation*} \label{rho_R}
\pi_0(\G)\stackrel{(\kappa,\cdots,\kappa)}{\longrightarrow}S_p\times\stackrel{q)}{\cdots}\times S_p\hookrightarrow S_{pq}.
\end{equation*}
\item[(ii)]
$\beta_{\Rb}:\pi_1(\G)\To(k^*)^{pq}_{\rho_{\Rb}}$ is the morphism
of $\pi_0(\G)$-modules defined by
$$
\beta_{\Rb}(u):=(\chi_1(g_1u),\ldots,\chi_1(g_pu),\ldots,\chi_q(g_1u),\ldots,\chi_q(g_pu)),\quad
u\in\pi_1(\G).
$$
\item[(iii)]
$c_{\Rb}:\pi_0(\G)\times\pi_0(\G)\To(k^*)^{pq}_{\rho_{\Rb}}$ is
the normalized 2-cochain defined by
\begin{align*}
\ \ \ \ c_{\Rb}&(g_i,g_j):= \\ &\ (\chi_1(\alpha(g_1,g_i,g_j)),\ldots,\chi_1(\alpha(g_p,g_i,g_j)),\ldots,\chi_q(\alpha(g_1,g_i,g_j)),\ldots,\chi_q(\alpha(g_p,g_i,g_j))
\end{align*}
for all $g_i,g_j\in\pi_0(\G)$.
\end{itemize}
\end{prop}

\begin{proof}
The proof is an easy but instructive exercise to become familiar
with the relationship between morphisms of 2-groups and the
associated triples described in \S~\ref{classes_iso_repr}. For
example, let us prove (i). We already know that for any $A\in{\rm
Obj}\Gg$ the functor $F_{\Rb}(A)$ basically amounts to permuting
the $\eta_{\chi,g}$, and we want to identify what this permutation
is. By definition we have
$$
F_{\Rb}(A)(\eta_{\chi,g'})(B)=\eta_{\chi,g'}(B\otimes A)=\left\{\begin{array}{ll} k, & \mbox{if $B\otimes A\in g'$} \\ 0, & \mbox{otherwise}.\end{array}\right.
$$
But $B\otimes A\in g'$ if and only if $B\in g'g^{-1}$, where
$g=[A]$. This means that $F_{\Rb}(A)(\eta_{\chi,g'})$ acts on
objects in exactly the same way as $\eta_{\chi,g'g^{-1}}$ and
consequently, we have
$$
F_{\Rb}(A)(\eta_{\chi,g'})\cong\eta_{\chi,g'g^{-1}}.
$$
Thus the morphism $\pi_0(\G)\To\pi_0(\ev_k^{\Gg})$ maps $g$ to the
isomorphism class of the permutation functor on $\ev_k^{\Gg}$
given by $\eta_{\chi,g'}\mapsto\eta_{\chi,g'g^{-1}}$, and under
our previous identification $\pi_0(\ev_k^{\Gg})\cong S_{pq}$ this
indeed corresponds to the morphism $\rho_{\Rb}$ defined above. We
leave to the reader the proof of (ii) and (iii). She/he can also
check that $\beta_{\Rb}$ indeed is a morphism of
$\pi_0(\G)$-modules and that $\partial c_{\Rb}=\beta_*(\alpha)$.
\end{proof}

\noindent In particular, although strictly speaking the regular
representation of a finite group $G$ is something different from
the regular representation of the associated discrete 2-group
$G[0]$, we see that the former is recovered as the equivalence class
of the later. For a finite one-object 2-group $A[1]$ we just get the set of
all characters of group $A$ as equivalence class of its regular representation.

Later on we shall use the triple
$(\rho_{\Rb},\beta_{\Rb},c_{\Rb})$ to get some more information on
the regular representation (see Example~\ref{exemple2} below).

\section{Categories of intertwiners.}

For any representations $\Fgg,\Fgg'$ let $\homs_{\G}(\Fgg,\Fgg')$,
or just ${\mathcal E}nd_{\G}(\Fgg)$ when $\Fgg=\Fgg'$, be the
associated category of intertwiners. It inherits an obvious
$k$-additive structure from the $k$-additive structures we have in
the underlying 2-vector spaces of each representation. In general,
however, it is not a 2-vector space because there may be no finite
basis of absolutely simple objects. For instance, ${\mathcal
E}nd_{\G}(\Ii)$ is equivalent to the category
$\rep_{\ev_k}(\pi_0(\G))$ of (finite dimensional) linear
representations of $\pi_0(\G)$ (see
Remark~\ref{remarca_repr_trivial} below). However, this is not
always a 2-vector space. Even if $\pi_0(\G)$ is finite, it may
lack to be a 2-vector space unless the field $k$ is algebraically
closed and of characteristic zero or prime to the order of
$\pi_0(\G)$.

At first sight, this is a little bit of a surprise when compared
to the corresponding situation for groups (finite or not), where
the set of intertwiners between any two finite dimensional linear
representations {\it always} is a finite dimensional vector space.
The difference arises from the fact that an intertwiner between
representations of a 2-group is not just a $k$-linear functor
between the underlying 2-vector spaces which satisfies some
additional conditions. That is to say, $\homs_{\G}(\Fgg,\Fgg')$ is
not a subcategory of $\homs_{k}(\Vv,\Vv')$. We further have the
all-important natural isomorphisms $\Phi(A)$ in
(\ref{condicio_intertwiner}) which come out as additional data we
are required to specify to completely define an intertwiner.

The purpose of this section is to prove that the same conditions
which ensure ${\mathcal E}nd_{\G}(\Ii)$ is a 2-vector space
(namely, $\pi_0(\G)$ finite and $k$ algebraically closed and of
characteristic zero or prime to the order of $\pi_0(\G)$) are
actually enough for the category $\homs_{\G}(\Fgg,\Fgg')$ to be a
2-vector space for any pair of representations $\Fgg$, $\Fgg'$. In
doing this we shall be able to describe explicitly a basis of
absolutely simple objects for these 2-vector spaces as well as a
method for computing the correspondings ranks out of the involved
representations. The proof is based on the geometric
interpretation of these categories of intertwiners given in
\cite{jE4} and recalled in \S~\ref{interpretacio_geometrica}.

All over this section various equivalences of categories are considered
whose explicit definitions will be needed in
Section~\ref{representabilitat}.

\subsection{The $k$-additive category $\homs_{\G}(\Fgg,\Fgg')$.}
\label{homs_G} Let $\Fgg=(\Vv,\F)$, $\Fgg'=(\Vv',\F')$. Then an
object in $\homs_{\G}(\Fgg,\Fgg')$ is given by a pair
$\xi=(H,\Phi)$ with $H:\Vv\To\Vv'$ a $k$-linear functor and
$\Phi=\{\Phi(A)\}_{A\in{\rm Obj}\Gg}$ a family of natural
isomorphisms of functors
\begin{equation} \label{2-isos_intertwiner}
\xymatrix{ \ar @{} [drr] |{\mbox{{\Large
        $\stackrel{\Phi(A)}{\Leftrightarrow}$}}}
\Vv\ar[d]_{F(A)}\ar[rr]^{H} & & \Vv'\ar[d]^{F'(A)} \\ \Vv\ar[rr]_{H} & &
\Vv'
}
\end{equation}
satisfying appropriate naturality and coherence conditions (see
\S~\ref{2-categoria_representacions}). In particular, if $R$ is
the matrix of ranks of $H$ (see \S~\ref{hom-categories_a_dev}),
the existence of such natural isomorphisms implies that $R$ is in
the obvious sense invariant under the action of $\pi_0(\G)$.

Among the objects in $\homs_{\G}(\Fgg,\Fgg')$ we have the {\it
zero intertwiner}, defined by the pair $(H_0,\Phi_0)$ with
$H_0:\Vv\To\Vv'$ ``the'' zero functor mapping all objects of $\Vv$
to a given zero object of $\Vv'$ and with all $\Phi_0(A)$ equal to
``identity'' natural transformations~\footnote{Strictly speaking,
the composites $H_0F(A)$ and $F'(A)H_0$ need not be equal. This is
the case if $F'(A)$ maps the given zero object of $\Vv'$ to
another zero object. Anyway, we always have a unique isomorphism
between both functors.}.

A morphism between two intertwiners $(H,\Phi)$ and
$(\tilde{H},\tilde{\Phi})$ is just a natural transformation
$\tau:H\Rightarrow\tilde{H}$ satisfying a naturality condition
which involves the 2-cells $\Phi(A)$ and $\tilde{\Phi}(A)$. It
follows that the zero intertwiner is a zero object of
$\homs_{\G}(\Fgg,\Fgg')$ and that $\homs_{\G}(\Fgg,\Fgg')$
inherits a $k$-linear structure from that existing in $\Vv'$ and
given by
$$
(\lambda\tau+\lambda'\tau')_V:=\lambda\tau_V+\lambda'\tau'_V,\qquad V\in{\rm Obj}\Vv,
$$
for any $\tau,\tau':H\Rightarrow\tilde{H}:\Vv\To\Vv'$ and any $\lambda,\lambda'\in k$. In particular, we have a forgetful $k$-linear functor
\begin{equation} \label{omega_F_F'}
\omega_{\Fgg,\Fgg'}:\homs_{\G}(\Fgg,\Fgg')\To\homs_{k}(\Vv,\Vv')
\end{equation}
mapping $(H,\Phi)$ to $H$ and equal to the identity on morphisms.
Notice, however, that this functor is neither injective nor
essentially surjective on objects and that it is a non-full
functor.

Biproducts in $\homs_{\G}(\Fgg,\Fgg')$ are obtained from the
biproducts in $\homs_k(\Vv,\Vv')$. More precisely, for objects
$(H,\Phi)$, $(\tilde{H},\tilde{\Phi})$ their biproduct is the pair
$(H\oplus\tilde{H},\Phi\oplus\tilde{\Phi})$ where
$H\oplus\tilde{H}$ is the biproduct in $\homs_k(\Vv,\Vv')$ (see
proof of Proposition~\ref{homs_interns}) and
$(\Phi\oplus\tilde{\Phi})(A)$ is given by the pasting
$$
\xymatrix{\ar @{} [drr] |{\mbox{{\Large
        $\cong$}}}
          \ar @{} [drrrrrr] |{\mbox{{\Large
        $\stackrel{\Phi(A)\oplus\tilde{\Phi}(A)}{\Leftrightarrow}$}}}
\Vv\ar@{=}[rr]\ar[d]_{F(A)} & & \Vv\ar@{=}[rr]\ar[d]_{h_r} & & \Vv\ar[rr]^{H\oplus\tilde{H}}\ar[d]^{h_l} & & \Vv'\ar[d]^{F'(A)}  \\ \Vv\ar[rr]_{H\oplus\tilde{H}} & & \Vv'\ar@{=}[rr] & & \Vv'\ar@{=}[rr] & & \Vv'\ar @{} [ull] |{\mbox{{\Large
        $\cong$}}}
}
$$
where $h_r:=(H\circ F(A))\oplus(\tilde{H}\circ F(A))$ and
$h_l:=(F'(A)\circ H)\oplus(F'(A)\circ\tilde{H})$. This makes sense
because composition of $k$-linear functors is $k$-bilinear and
hence, distributes over biproducts in a canonical way. We leave to
the reader checking that the pair $(H\oplus H',\Phi\oplus\Phi')$
so defined is indeed a new intertwiner between $\Fgg$ and $\Fgg'$.

\subsection{Notation.}
If $\Fgg_1\simeq\Fgg'_1$ and $\Fgg_2\simeq\Fgg'_2$ we clearly have
$\homs_{\G}(\Fgg_1,\Fgg'_1)\simeq\homs_{\G}(\Fgg_2,\Fgg'_2)$. To
emphasize this, in the rest of this section we denote the
intertwining hom-categories by
$$
\homs_{\G}(\Fgg,\Fgg')\ \equiv\ \Hh\left(\begin{array}{c} n,\rho,\beta,c \\ n',\rho',\beta',c'\end{array}\right),
$$
or just $\Hh(n,\rho,\beta,c)$ when both representations are the
same (up to equivalence). The reader may think of these categories
$\Hh\left(\begin{array}{c} n,\rho,\beta,c \\
n',\rho',\beta',c'\end{array}\right)$ as the hom-categories
between specific representatives we have fixed once and for all
for each equivalence class of representations. For instance, those
described in \S~\ref{classes_iso_repr}.

\subsection{Geometric description of the categories of intertwiners.}
\label{interpretacio_geometrica}

Let $G$ be any group and $X$ a right $G$-set. We shall denote by
$F(X,k^*)$ the (multiplicative) abelian group of all $k^*$-valued
functions on $X$. When we speak of 2-cocycles of $G$ with values
in $F(X,k^*)$ we always assume $F(X,k^*)$ to be equipped with the
$G$-module structure
$$
(g\cdot f)(x)=f(xg),\qquad x\in X.
$$
Let $z$ be a normalized 2-cocycle of $G$ with values
in $F(X,k^*)$ (i.e. a 2-cocycle such that $z(g,e)=z(e,g)=1$ for any
$g\in G$, where $1$ denotes the unit of $F(X,k^*)$).

Given $(G,X,z)$ as above, we denote by $\Vv ect_{G,z}(X)$ the
corresponding category of $z$-projective $G$-equivariant vector
bundles over $X$. Objects are given by triples $(E,p,\Theta)$
with $(E,p)$ a finite rank vector bundle $p:E\To X$ over $X$, and
$\Theta:E\times G\To E$ a $z$-projective right $G$-action making
$p$ a $G$-equivariant map and whose restriction to fibers is
$k$-linear. Thus if we denote by $\theta(x,g):E_x\To E_{xg}$ the
$k$-linear isomorphisms defined by the restriction of $\Theta$ to
$E_x\times\{g\}$ we have
\begin{align} \label{z-projectiva_1}
\theta(x,gg')&=z(g,g')(x)\ \theta(xg,g')\circ\theta(x,g) \\  \label{z-projectiva_2} \theta(x,e)&={\id}_{E_x}
\end{align}
for all $g,g'\in G$ and $x\in X$. A morphism
$\phi:(E,p,\Theta)\To(E',p',\Theta')$ between two such triples is
an action preserving morphism of vector bundles, hence a family
$$
\phi=\{ \phi_x:E_x\To E'_x \}_{x\in X}
$$
of $k$-linear maps such
that
\begin{equation} \label{morfisme_fibrats}
\phi_{xg}\circ \theta(x,g)=\theta'(x,g)\circ\phi_x
\end{equation}
for all $g\in G$ and $x\in X$. Composition is the obvious one.

Observe that in writting $\Vv ect_{G,z}(X)$ we do not make
explicit the field $k$. But it is there. Actually, $\Vv
ect_{G,z}(X)$ is a $k$-additive category. The $k$-linear structure
is the obvious one,
the zero vector bundle equipped with its unique $z$-projective
right $G$-action is a zero object, and
$(E,p,\Theta)\oplus(E',p',\Theta')$ is the usual direct sum of
vector bundles equipped with the $z$-projective action
$$
(\theta\oplus\theta')(x,g):E_x\oplus E'_x\To E_{xg}\oplus
E'_{xg}
$$
defined by
$$
(\theta\oplus\theta')(x,g)(v_x+v'_x):=\theta(x,g)(v_x)+\theta'(x,g)(v'_x),\quad v_x\in E_x,\ v'_x\in E'_x.
$$
As we will see later, it is even a 2-vector space under suitable
assumptions.

Let now $(n,\rho,\beta,c)$ and $(n',\rho',\beta',c')$ be quadruples
of the kind described in \S~\ref{classes_iso_repr}. The group
morphisms $\rho$ and $\rho'$ induce a right action of $\pi_0(\G)$
on $X(n',n):=[n']\times[n]$ given by
$$
(i',i)\cdot g=(\rho'(g^{-1})(i'),\rho(g^{-1})(i)),\quad g\in G.
$$
Let us denote by $\Lambda(n,\rho,\beta;n',\rho',\beta')$ the
corresponding set of {\it intertwining} $\pi_0(\G)$-orbits, i.e.
orbits $X_{\lambda}$ such that $\beta_i=\beta'_{i'}$ for
all~\footnote{Actually, it is easy to see that this condition
holds for all points in $X_{\lambda}$ if it holds for some
(arbitrary) point $(i',i)\in X_{\lambda}$.} $(i',i)\in
X_{\lambda}$. Finally, for each intertwining $\pi_0(\G)$-orbit
$X_{\lambda}$ a normalized 2-cocycle $z_{\lambda}\in
Z^2(\pi_0(\G),F(X_{\lambda},k^*))$ is defined by
$$
z_{\lambda}(g_1,g_2)(i',i)=\frac{c'(g_1,g_2)_{i'}}{c(g_1,g_2)_i}
$$
for all $g_1,g_2\in\pi_0(\G)$ and $(i',i)\in X_{\lambda}$. Then we have the
following. 

\begin{thm}[\cite{jE4}] \label{hom-categories}
There is an equivalence of $k$-additive categories
\begin{equation} \label{hom_categories}
\Hh\left(\begin{array}{c} n,\rho,\beta,c \\ n',\rho',\beta',c'\end{array}\right)\simeq
\prod_{X_{\lambda}\in\Lambda(n,\rho,\beta;n',\rho',\beta')} \Vv ect_{\pi_0(\G),z_{\lambda}}(X_{\lambda}).
\end{equation}
\end{thm}
\noindent For later use, let us recall from \cite{jE4} how this equivalence
works. Let $(H,\Phi)$ be any intertwiner, and let $R=(r_{i'i})$ be the
matrix of ranks of the functor $H$. As mentioned before, $R$ is
invariant under the action of $\pi_0(\G)$. Hence associated to
each orbit $X_{\lambda}\subset X(n',n)$ we have a well defined
nonnegative integer $d_{\lambda}$ (the common value of the
corresponding entries in $R$). This gives the rank of the vector
bundle $(E(\lambda),p(\lambda))$ over $X_{\lambda}$, and it is
easy to see that this rank is necessarily zero unless
$X_{\lambda}$ is an intertwining orbit. Let us assume without loss
of generality that
$$
E(\lambda)=\coprod_{x\in X_{\lambda}}k^{d_{\lambda}}
$$
and that $p(\lambda)$ is the obvious projection. The
$z_{\lambda}$-projective action $\Theta(\lambda)$ is now
determined by the natural isomorphisms $\Phi(A)$. To be explicit,
let us think of the left hand side of (\ref{hom_categories}) as
the category of intertwiners between the representations
$\Fgg(n,\rho,\beta,c),\Fgg'(n',\rho',\beta',c')$ described in
\S~\ref{classes_iso_repr}. Thus the underlying 2-vector spaces
$\Vv,\Vv'$ are of the form $\ev_k^r$ in both representations, and
$F(A)$ and $F'(A)$ are permutation functors $P_{\rho[A]}$ and
$P_{\rho'[A]}$, respectively. In this case $\Phi(A)$ is a natural
isomorphism
$$
\Phi(A):P_{\rho'[A]}H\Rightarrow HP_{\rho[A]}:\ev_k^n\To\ev_k^{n'}
$$
and hence, it is given by an $n'\times n$ matrix whose
$(i',i)^{th}$-entry is itself an invertible matrix
$$
\Phi(A)_{i'i}\in{\rm GL}(r_{i',\rho[A](i)},k)
$$
with entries in $k$ if $r_{i',\rho[A](i)}\neq 0$ (otherwise, it is
the empty matrix; see \S~\ref{hom-categories_a_dev}). Then the
linear isomorphisms $\theta((i',i),g):k^{d_{\lambda}}\To
k^{d_{\lambda}}$ defining the action $\Theta(\lambda)$ are those
which in canonical bases are given by the invertible matrices
$$
\theta((i',i),g)\ \stackrel{\mbox{\small canonical
bases}}{\longleftrightarrow}\
\left(\Phi(A)_{i',\rho(g^{-1})(i)}\right)^{-1}\in{\rm
GL}(d_{\lambda},k),\quad (i',i)\in X_{\lambda}
$$
for any $A$ such that $[A]=g$. Then (\ref{hom_categories}) maps the object
$(H,\Phi)$ to
$$
(E(\lambda),p(\lambda),\Theta(\lambda))_{X_{\lambda}}\ \in\ {\rm
Obj}\left(\prod_{X_{\lambda}\in\Lambda(n,\rho,\beta;n',\rho',\beta')}
\Vv ect_{\pi_0(\G),z_{\lambda}}(X_{\lambda})\right).
$$
The action on morphisms is as follows. Let
$\tau:H\Rightarrow\tilde{H}:\Vv\To\Vv'$ be a morphism from
$(H,\Phi)$ to $(\tilde{H},\tilde{\Phi})$ for any intertwiners
$(H,\Phi),(\tilde{H},\tilde{\Phi}):(\Vv,\F)\To(\Vv',\F')$. As
pointed out before, $\tau$ is completely given by its components
$\tau_{V_i}:H(V_i)\To\tilde{H}(V_i)$ on a basis
$\{V_1,\ldots,V_n\}$ of $\Vv$, and each of these components is in
turn described by $n'$ matrices $M_{i'i}\in{\rm
Mat}_{\tilde{r}_{i'i}\times r_{i'i}}(k)$, $i'=1,\ldots,n'$ (cf.
proof of Proposition~\ref{homs_interns}). Then $\tau$ gets mapped
to the morphism $\phi=(\phi(\lambda))_{\lambda}$ whose
$X_{\lambda}$-component
$$
\phi(\lambda):(E(\lambda),p(\lambda),\Theta(\lambda))\longrightarrow(\tilde{E}(\lambda),\tilde{p}(\lambda),\tilde{\Theta}(\lambda))
$$
is the morphism in $\Vv ect_{\pi_0(\G),z_{\lambda}}(X_{\lambda})$
given on fibers by these matrices $M_{i'i}$. More precisely, if
$(i',i)\in X_{\lambda}$ the map
$$
\phi(\lambda)_{(i',i)}:E(\lambda)_{(i',i)}=k^{d_{\lambda}}\longrightarrow
k^{\tilde{d}_{\lambda}}= \tilde{E}(\lambda)_{(i',i)}
$$
is the $k$-linear map given in canonical bases by the matrix $M_{i'i}$. The morphism 
$\phi$ so defined satisfies (\ref{morfisme_fibrats}) because of
the above mentioned condition on $\tau$ involving the 2-cells
$\Phi(A)$ and $\tilde{\Phi}(A)$ and ensuring that $\tau$ is indeed
a 2-intertwiner between $(H,\Phi)$ and $(\tilde{H},\tilde{\Phi})$
(recall that the functor (\ref{omega_F_F'}) is non-full!).

\begin{rem}{\rm
In \cite{jE4} we proved that this functor is an equivalence of
categories. In fact the functor is $k$-linear and hence, the
equivalence is of $k$-additive categories. Indeed, any $k$-linear
functor between $k$-additive categories automatically preserves
biproducts; see \cite{sM98}, p. 197 where this is shown for the
case the commutative ring $k$ is $\mathbb{Z}$. }
\end{rem}

\begin{rem} \label{remarca_repr_trivial}{\rm
If $n=n'=1$, and $\beta=\beta'$ and $c=c'$ are the trivial maps
$\pi_0(\G)^3\To k^*$ and $\pi_0(\G)^2\To k^*$, respectively, we
have $\Fgg,\Fgg'\simeq\Ii$. In this case, the right hand side of
(\ref{hom_categories}) indeed reduces to the category
$\rep_{\ev_k}(\pi_0(\G))$. In fact, the equivalence is in this
case as monoidal categories when $\mathcal{E}nd_{\G}(\Ii)$ comes
equipped with the monoidal structure induced by the composition of
endomorphisms and $\rep_{\ev_k}(\pi_0(\G))$ with the usual tensor
product of representations. This implies that we shall have no
analog of Schur's lemma, at least in its usual version. Indeed,
whatever definition we adopt for the irreducible representations
in this 2-category setting, the trivial representation $\Ii$
should be such a representation. But linear representations of
groups, in our case of $\pi_0(\G)$, have no inverse with respect
to tensor product. Therefore $\Ii$ will be an irreducible
representation with lots of non-invertible nonzero endomorphisms.
}
\end{rem}

\subsection{The categories $\Vv ect_{G,z}(X)$ for a transitive $G$-set $X$.}

It readily follows from Theorem~\ref{hom-categories} and
Proposition~\ref{homs_interns} that $\Hh\left(\begin{array}{c}
n,\rho,\beta,c \\ n',\rho',\beta',c'\end{array}\right)$ will be a
2-vector space when all $k$-additive categories $\Vv
ect_{\pi_0(\G),z_{\lambda}}(X_{\lambda})$ are 2-vector spaces. To
prove that these categories are indeed 2-vector spaces we shall
take advantage of the fact that all $\pi_0(\G)$-sets $X_{\lambda}$ are
transitive to get a more
elementary description of them.

Let us start with the following observation.
\begin{lem}
Let $X$ be a {\it transitive} (right) $G$-set and for any $x\in X$
let $G_x\subset G$ be the stabilizer of $x$. Then any 2-cocycle
$z\in Z^2(G,F(X,k^*))$ gives rise to 2-cocycles $z_x,\hat{z}_x\in
Z^2(G_x,k^*)$ defined by
\begin{eqnarray*}
z_x(g_1,g_2)&:=&z(g_1,g_2)(x)
\\ \hat{z}_x(g_1,g_2)&:=&z_(g^{-1}_2,g^{-1}_1)(x)
\end{eqnarray*}
for any $g_1,g_2\in G_x$. Here $k^*$ is assumed to be equipped
with the trivial $G_x$-module structure. Furthermore, $z_x$ and
$\hat{z}_x$ are normalized when $z$ is normalized.
\end{lem}
\begin{proof}
An easy computation shows that
\begin{align*}
\partial z_x(g_1,g_2,g_3)&=\partial z(g_1,g_2,g_3)(x) \\
\partial \hat{z}_x(g_1,g_2,g_3)&=\left(\partial z(g^{-1}_3,g_2^{-1},g_1^{-1})\right)^{-1}(x)
\end{align*}
for all $g_1,g_2,g_3\in G_x$.
\end{proof}
\noindent Recall that for any group $H$ and any normalized
2-cocycle $z\in Z^2(H,k^*)$ a $z$-{\it projective representation}
of $H$ (or projective representation with {\it central charge}
$z$) is a vector space $V$ together with a map $\psi:H\To GL(V)$
such that $\psi(e)={\id}_V$ and
$$
\psi(h_1h_2)=z(h_1,h_2)\ \psi(h_1)\circ\psi(h_2)
$$
for all $h_1,h_2\in H$. These representations are the objects of a
category $\rep_{z}(H)$ whose morphisms~\footnote{Let us remark
that there exists a more general notion of morphism between
projective representations (with the same or with different
central charges) called {\it projective morphisms}. These are
given by a $k$-linear map $f:V\To V'$ together with a map
$\mu:H\To k^*$ such that $\mu(1)=1$ and
$f\circ\psi(h)=\mu(h)\psi'(h)\circ f$. When $z=z'$ it follows that
$\mu$ is a homomorphism. Clearly we have embeddings
$\rep_z(H)\hookrightarrow\mathcal{P}\rep_z(H)\hookrightarrow\mathcal{P}\rep(H)$,
where $\mathcal{P}\rep_z(H)$ denotes the category of
$z$-projective representations of $H$ with the projective
morphisms, and $\mathcal{P}\rep(H)$ the category of $z$-projective
representations of $H$ for arbitrary 2-cocycles $z$, and
projective morphisms between them.} are $k$-linear maps $f:V\To
V'$ such that $f\circ\psi(h)=\psi'(h)\circ f$ for all $h\in H$. In
particular, when $z$ is trivial we recover the category of linear
representations of $H$.

\begin{prop} \label{prop1}
Let $G$ be a group, $X$ a transitive (right) $G$-set and $z\in
Z^2(G,F(X,k^*))$ a normalized 2-cocycle. Then for any $x_0\in X$
we have an equivalence of $k$-additive categories
$$
\Vv ect_{G,z}(X)\simeq\rep_{\hat{z}_{x_0}}(G_{x_0}),
$$
where $G_{x_0}\subset G$ is the stabilizer (or isotropy subgroup) of $x_0$.
\end{prop}
\begin{proof}
For any object $(E,p,\Theta)$ of $\Vv ect_{G,z}(X)$ let
$\psi:G_{x_0}\To {\rm GL}(E_{x_0})$ be defined by
$$
\psi(g):=\theta(x_0,g^{-1}),\qquad g\in G_{x_0}.
$$
It readily follows from (\ref{z-projectiva_1}) and
(\ref{z-projectiva_2}) that $\psi$ is a $z_{x_0}$-projective
representation of $G_{x_0}$. Moreover, from
(\ref{morfisme_fibrats}) it follows that the $x_0$-component
$\phi_{x_0}:E_{x_0}\To E'_{x_0}$ of any morphism
$\phi:(E,p,\Theta)\To(E',p',\Theta')$ in $\Vv ect_{G,z}(X)$ is an
intertwiner between the corresponding representations $\psi$ and
$\psi'$. This defines a $k$-linear functor
$$
F:\Vv ect_{G,z}(X)\To\rep_{\hat{z}_{x_0}}(G_{x_0}),
$$
and we claim that this functor is an equivalence of categories.

Indeed, transitivity of $X$ together with (\ref{morfisme_fibrats})
show that any morphism $\phi$ in $\Vv ect_{G,z}(X)$ is uniquely
determined by its $x_0$-component $\phi_{x_0}$ and moreover, that any
intertwiner $f:E_{x_0}\To E'_{x_0}$ between $\psi$ and $\psi'$ is
the $x_0$-component of such a $\phi$ (i.e. it can be extended to a
whole morphism $\phi$ between $(E,p,\Theta)$ and
$(E',p',\Theta')$). Hence $F$ is fully faithful.

To prove $F$ is essentially surjective, let $\psi:G_{x_0}\To {\rm
GL}(V)$ be any $\hat{z}_{x_0}$-projective representation. An
object of $\Vv ect_{G,z}(X)$ can be built out of it as follows.
Let us fix representatives $\Rr=\{g_1,\ldots,g_r\}$ of the right
cosets of $G_{x_0}/G$, with $g_1=e$ as representative of
$G_{x_0}$. Set $E:=\coprod_{x\in X} V$ and let $p:E\To X$ be the
obvious projection. Because of the transitivity of $X$, there
exist unique $g_i,g_{i'}\in\Rr$ and $\tilde{g}\in G_{x_0}$ such
that
\begin{equation} \label{gi_tildeg_gi'}
x=x_0g_i,\qquad g_ig=\tilde{g}g_{i'}.
\end{equation}
Then for any pair $(x,g)\in X\times G$ let $\theta(x,g):E_{x}\To
E_{xg}$ be the $k$-linear isomorphism defined by
\begin{equation} \label{definicio_theta}
\theta(x,g):=\frac{z(\tilde{g},g_{i'})(x_0)}{z(g_i,g)(x_0)}\psi(\tilde{g}^{-1}).
\end{equation}
Let us see that the pair $(E,p)$ together with these maps indeed
define an object of $\Vv ect_{G,z}(X)$. If $g=e$ we have
$\tilde{g}=e$ and $g_{i'}=g_i$. Hence (\ref{z-projectiva_2}) holds
because $z$ is normalized. To prove (\ref{z-projectiva_1}) let
$g_j,g_{j'}\in\Rr$ and $\hat{g}\in G_{x_0}$ be uniquely defined by
\begin{equation} \label{gj_hatg_gj'}
xg=x_0g_j,\qquad g_jg'=\hat{g}g_{j'}.
\end{equation}
Hence
$$
\theta(xg,g')=\frac{z(\hat{g},g_{j'})(x_0)}{z(g_j,g')(x_0)}\psi(\hat{g}^{-1}).
$$
Similarly, let $g_{i''}\in\Rr$ and $\overline{g}\in G_{x_0}$ be uniquely defined by
\begin{equation} \label{gi''_overlineg}
g_igg'=\overline{g}g_{i''}
\end{equation}
so that the left  hand side of (\ref{z-projectiva_1}) is
$$
\theta(x,gg')=\frac{z(\overline{g},g_{i''})(x_0)}{z(g_i,gg')(x_0)}\psi(\overline{g}^{-1}).
$$
Thus we have to prove that
\begin{equation} \label{a_provar}
\frac{z(\overline{g},g_{i''})(x_0)}{z(g_i,gg')(x_0)}\psi(\overline{g}^{-1})=z(g,g')(x_0g_i)\frac{z(\hat{g},g_{j'})(x_0)}{z(g_j,g')(x_0)}\frac{z(\tilde{g},g_{i'})(x_0)}{z(g_i,g)(x_0)}\psi(\hat{g}^{-1})\psi(\tilde{g}^{-1})
\end{equation}
To show this, note first that not all of elements
$g_i,g_{i'},g_{i''},g_j,g_{j'}\in\Rr$ are independent, and the
same is true for the elements $\tilde{g},\hat{g},\overline{g}\in
G_{x_0}$. Thus from (\ref{gi_tildeg_gi'}) and (\ref{gj_hatg_gj'})
we have
$$
x_0g_j=xg=x_0g_ig=x_0\tilde{g}g_{i'}=x_0g_{i'}
$$
so that $g_j=g_{i'}$. Using now (\ref{gi''_overlineg}) it follows that
$$
\overline{g}g_{i''}=g_igg'=\tilde{g}g_{i'}g'=\tilde{g}g_jg'=\tilde{g}\hat{g}g_{j'}
$$
so that $\overline{g}=\tilde{g}\hat{g}$ and $g_{i''}=g_{j'}$. Moreover we have
$$
\psi(\tilde{g}^{-1})\psi(\hat{g}^{-1})=\frac{1}{z(\hat{g},\tilde{g})(x_0)}\ \psi(\tilde{g}^{-1}\hat{g}^{-1})
$$
because $\psi$ is $\hat{z}_{x_0}$-projective. Putting all these
facts together we see that (\ref{a_provar}) reduces to
\begin{equation}\label{reescriptura}
z(\tilde{g}\hat{g},g_{i''})z(g_{i'},g')z(g_i,g)z(\tilde{g},\hat{g})_{|_{x_0}}=(g_iz(g,g'))z(\hat{g},g_{i''})z(\tilde{g},g_{i'})z(g_i,gg')_{|_{x_0}},
\end{equation}
where we have used that $z(g,g')(x_0g_i)=(g_iz(g,g'))(x_0)$. Now
by the 2-cocycle condition on $z$ we have
\begin{align*}
z(\tilde{g}\hat{g},g_{i''})(x_0)&=z(\hat{g},g_{i''})z(\tilde{g},\hat{g}g_{i''})z(\tilde{g},\hat{g})^{-1}\ _{|_{x_0}} \\
z(g_i,gg')(x_0)&=(g_iz(g,g'))^{-1}z(g_ig,g')z(g_i,g)_{|_{x_0}}.
\end{align*}
In the first equality we have used that $\tilde{g}\in G_{x_0}$ so
that $(\tilde{g}z(\hat{g},g_{i''}))(x_0)=z(\hat{g},g_{i''})(x_0)$.
Putting this into (\ref{reescriptura}) and using that
$g_ig=\tilde{g}g_{i'}$ and $\hat{g}g_{i''}=g_{i'}g'$ shows that
(\ref{reescriptura}) holds. To finish the proof it remains to see
that the object $(E,p,\Theta)$ of $\Vv ect_{G,z}(X)$ we have
constructed in this way out of $\psi$ indeed gets mapped by the
functor $F$ to a $\hat{z}_{x_0}$-projective representation
equivalent to $\psi$. In fact, it gets mapped to $\psi$ because
for any $g\in G_{x_0}$ we have
$$
F(E,p,\Theta)(g)=\theta(x_0,g^{-1})\stackrel{(\ref{definicio_theta})}{=}\frac{z(g^{-1},e)(x_0)}{z(e,g^{-1})(x_0)}\psi(g)=\psi(g).
$$
Here we use that the $g_i,g_{i'},\tilde{g}$ in
(\ref{gi_tildeg_gi'}) are in this case given by $g_i=g_{i'}=e$ and
$\tilde{g}=g^{-1}$ because $x=x_0$ and $g\in G_{x_0}$.
\end{proof}

\begin{rem}{\rm
We have shown that $F$ is surjective on objects, not just
essentially surjective. However, $F$ is not an isomorphism of
categories because it is not injective on objects. Indeed, to
construct a preimage of $\psi$ we need to choose representatives
for the right cosets in $G_{x_0}/G$, and different choices will
give isomorphic, but not equal, objects in $\Vv ect_{G,z}(X)$
which get mapped to $\psi$ by the functor $F$. Note also that any
pseudoinverse of $F$ will map an intertwiner $f:V\To V'$ in
$\rep_{z_{x_0}}(G_{x_0})$ to the unique morphism
$\phi:\coprod_{x\in X}V\To\coprod_{x\in X}V'$ whose restriction to
the fiber over $x_0$ is $f$. }
\end{rem}
\noindent The following is an immediate consequence of the
previous result and the obvious fact that $\Rr ep_z(1)=\ev_k$.
\begin{cor} \label{cas_particular_torsors}
Let $X$ be a $G$-torsor (i.e. a transitive $G$-set with trivial
stabilizers). Then we have an equivalence of $k$-additive
categories
$$
\Vv ect_{G,z}(X)\simeq\ev_k
$$
for any normalized 2-cocycle $z\in Z^2(G,F(X,k^*))$.
\end{cor}
\noindent In particular, we conclude that when $X$ is a $G$-torsor
the isomorphism class of any object $(E,p,\Theta)$ of $\Vv
ect_{G,z}(X)$ is completely given by its rank $d\geq 0$. A
specific representative in this isomorphism class is the triple
$(E(d),p(d),\Theta(d))$ with $E(d)=\coprod_{x\in X}k^d$, $p(d)$
the obvious projection, and $\Theta(d)$ given by
\begin{equation} \label{theta_cas_torsor}
\theta(d)(x,g)=z(\overline{g},g)(x_0)^{-1}\ {\id}_{k^d}
\end{equation}
for any $x_0\in X$ and $\overline{g}\in G$ the unique such that
$x=x_0\overline{g}$; cf. (\ref{definicio_theta}). Moreover, a
$k$-linear map $f:k^d\To k^{\tilde{d}}$ corresponds to the
morphism
$\phi(f):(E(d),p(d),\Theta(d))\To(E(\tilde{d}),p(\tilde{d}),\Theta(\tilde{d}))$
whose components $\phi(f)_x$ are all equal to $f$. Indeed, as
pointed out in the above remark we have $\phi(f)_{x_0}=f$, while
the other components follow from (\ref{theta_cas_torsor}) and
(\ref{morfisme_fibrats}).

\subsection{Review on projective representations and modules over arbitrary semisimple algebras} \label{representacions_projectives}
In this subsection we recall a few well known facts from the
theory of projective representations of finite groups and more
generally, of modules over a semisimple $k$-algebra. The aim is to
see that, under appropriate assumptions on the field $k$, the
corresponding categories (for a given central charge in the case
of projective representations) are 2-vector spaces, and to explain
how their ranks can be computed. This result generalizes
Example~\ref{ex2} in \S~\ref{definicio_2ev} and allows us to prove
that all hom-categories in $\repg_{\dev_k}(\G)$ indeed are
2-vector spaces under the appropriate assumptions. We refer the
reader to \cite{gK85} for the theory of projective representations
of a finite group and to \cite{sW03} for the general case.

Ordinary linear representations of a finite
group $G$ are the same as (left) modules over the group algebra
$k[G]$ and moreover, $k[G]$ is a semisimple $k$-algebra when $k$ is
algebraically closed of characteristic zero or prime to the order
of $G$. These are the two basic facts which prove that the
category $\rep_{\ev_k}(G)$ of Example~\ref{ex2} is a 2-vector
space.

More generally, let $A$ be any finite dimensional semisimple
$k$-algebra, with $k$ algebraically closed. Then each finite
dimensional $A$-module decomposes as a finite direct sum of
irreducible $A$-modules, and this decomposition is unique up to
isomorphism and permutation of the factors (see Theorem~2.2 in
Chapter~2 of \cite{sW03}). Moreover, irreducible modules are
absolutely simple in our sense above, and there are only finitely
many isomorphism classes of them (Lemma~2.1 and Corollary~2.15 in
Chapter 2 of {\it loc. cit}). Briefly, the category $A$-$\modb$ of
$A$-modules is a 2-vector space with basis of absolutely simple
objects any set of representatives of the irreducible modules. If
$A=k[G]$ the condition on $k$ to be of characteristic zero or
prime to the order of $G$ is just the necessary and suficient
condition for $k[G]$ to be semisimple (this is the famous
Maschke's theorem; see Theorem~1.14 in Chapter~3 of \cite{sW03}).

Let us now consider $z$-projective representations for a given
normalized 2-cocycle $z$. The first remark is that these
representations are the same as modules over the {\it twisted
group algebra} $k[G]_z$. This is the $k$-algebra with the same
underlying space as $k[G]$ but with multiplication given by
$$
e_{g}e_{g'}:=z(g,g')e_{gg'},\qquad g,g'\in G
$$
(cf. Chapter~3, \S~2 of \cite{gK85}). The second remark is that
twisted group algebras are also semisimple $k$-algebras when $k$
is of characteristic zero or prime to the order of $G$. The proof
is essentially the same as for $k[G]$ (see Theorem~2.10 in
Chapter~3 of \cite{gK85}). Therefore, always under the assumption
that $k$ is algebraically closed, $\rep_z(G)$ is a 2-vector space
with basis of absolutely simple objects any set of representatives
of the irreducible modules.

What about ranks? Let $A$ be an arbitrary finite dimensional
semisimple $k$-algebra, and let $\{M_i,i\in I\}$ be any set of
representatives of the isomorphism classes of irreducible
$A$-modules. Then it is shown that $|I|={\rm dim}_kZ(A)$, where
$Z(A)$ denotes the center of $A$, and that $A\cong\oplus_{i\in
I}n_iM_i$, with the $n_i\geq 0$ such that
$$
n_i={\rm dim}_kM_i
$$
and
$$
{\rm dim}_kA=\sum_{i\in I}n_i^2
$$
(cf. Corollary~2.24 in Chapter~2 of \cite{sW03}). In particular, the rank of $A$-$\modb$ as a 2-vector space is equal to the dimension over $k$ of the center of $A$.

This reduces the problem of computing the rank of the 2-vector
space $\rep_z(G)$ to that of computing the dimension over $k$ of
the center of $k[G]_z$. If $z=1$ we recover the usual group
algebra $k[G]$, and it is well known that a $k$-basis of its
center is given by the elements $c_i=\sum_{g\in C_i}e_g$,
$i=1,\ldots,t$, if $C_1,\ldots,C_t$ are the conjugacy classes of
$G$ (Lemma~3.2 in Chapter~3 of \cite{sW03}). This gives the
statement in Example~\ref{ex2}. The answer for an arbitrary
nontrivial normalized 2-cocycle $z$ can be found in \cite{gK85}
(Chapter~3, \S~6). If $k$ is of characteristic zero the answer is
the following. An element $g\in G$ is called $z$-{\it regular} if
$e_{g'}e_g=e_ge_{g'}$ for all $g'\in C_g(G)$, the centralizer of
$g$ in $G$. The product here takes place in the twisted group
algebra $k[G]_z$. In other words, $g\in G$ is $z$-regular iff
$z(g,g')=z(g',g)$ for all $g'\in C_g(G)$. It is easy to check that
if $g\in G$ is $z$-regular then so is any conjugate of $g$. Hence
it makes sense to speak of the $z$-regular conjugacy classes of
$G$. Let $C_1,\ldots,C_{t'}$ be all the $z$-regular conjugacy
classes of $G$. Then it is shown that the elements $c_i=\sum_{g\in
C_i}e_g$, $i=1,\ldots,t'$, constitute a $k$-basis of $Z(k[G]_z)$.

Therefore we have the following generalization of
Example~\ref{ex2}, and a restatement of the above mentioned
results from \cite{gK85}.

\begin{prop} \label{prop2}
Let $k$ be an algebraically closed field of characteristic zero,
$G$ a finite group and $z\in Z^2(G,k^*)$ a normalized 2-cocycle.
Then $\rep_{z}(G)$ is a 2-vector space of rank the number of
$z$-regular conjugacy classes of $G$, a basis of absolutely simple
objects being given by any set of representatives of the
isomorphism classes of irreducible $z$-projective representations.
\end{prop}

\subsection{Main theorem of the section} \label{teorema_principal}
Next result readily follows now from Theorem~\ref{hom-categories} and Propositions~\ref{prop1} and \ref{prop2}.

\begin{thm}
Let $\G$ be any 2-group and
$(n,\rho,\beta,c),(n',\rho',\beta',c')$ quadruples of the kind
described in \S~\ref{classes_iso_repr}. Then there is a $k$-linear
equivalence of categories
$$
\Hh\left(\begin{array}{c} n,\rho,\beta,c \\ n',\rho',\beta',c'\end{array}\right)\simeq\prod_{X_{\lambda}\in\Lambda(n,\rho,\beta;n',\rho',\beta')} \mathcal{R}ep_{\hat{z}_{\lambda}}(G_{i'_{\lambda}i_{\lambda}}),
$$
where $(i'_{\lambda},i_{\lambda})$ is any point of $X_{\lambda}$,
$G_{i'_{\lambda}i_{\lambda}}\subseteq\pi_0(\G)$ the corresponding
stabilizer, and $\hat{z}_{\lambda}\in
Z^2(G_{i'_{\lambda}i_{\lambda}},k^*)$ the normalized 2-cocycle
defined by
$$
\hat{z}_{\lambda}(g_1,g_2):=\frac{c'(g_2^{-1},g_1^{-1})_{i'_{\lambda}}}{c(g_2^{-1},g_1^{-1})_{i_{\lambda}}}.
$$
Moreover, when $\G$ is essentially finite and $k$ is algebraically
closed and of characteristic zero or prime to the order of
$\pi_0(\G)$ these $k$-linear categories are 2-vector spaces.
\end{thm}
Note that the finiteness of $\pi_0(\G)$ ensures that each
$k$-additive category
$\mathcal{R}ep_{\hat{z}_{\lambda}}(G_{i'_{\lambda}i_{\lambda}})$
is a 2-vector space, while that of $\pi_1(\G)$ ensures that there
is a finite number of intertwining orbits in $X(n',n)$ and hence,
a finite number of terms in the above product.

The following special case is important for what follows. In
particular, it is used in \S~\ref{seccio_functor_universal} to
prove the representability of the forgetful 2-functor by the
regular representation (cf. also
Corollary~\ref{cas_particular_torsors}).

\begin{cor}\label{corol.lari}
Let $\G$ be essentially finite and $k$ algebraically closed and of
characteristic zero or prime to the order of $\pi_0(\G)$. In the
above notations, let us assume that
$(n,\rho,\beta,c),(n',\rho',\beta',c')$ are such that all
intertwining orbits $X_{\lambda}$ are $\pi_0(\G)$-torsors. Then we
have a $k$-linear equivalence of categories
\begin{equation} \label{cas_torsors}
\Hh\left(\begin{array}{c} n,\rho,\beta,c \\ n',\rho',\beta',c'\end{array}\right)\simeq\ev_k^N,
\end{equation}
where $N$ is the number of intertwining orbits in $X(n',n)$.
\end{cor}

\begin{rem}\label{remarca_matriu_rangs}{\rm
From the whole discussion above it follows that the equivalence
(\ref{cas_torsors}) goes as follows (from right to left). On the
one hand, an object $(k^{d_1},\ldots,k^{d_N})$ of $\ev_k^N$ is
mapped to any intertwiner $(H,\Phi):(\Vv,\F)\To(\Vv',\F')$ whose
functor $H$ has a matrix of ranks $R=(r_{i'i})$ given by
$r_{i'i}=d_{\lambda}$ for all $(i',i)\in X_{\lambda}$ (this
completely determines the intertwiner up to isomorphism). In
particular, a basis of
$\Hh\left(\begin{array}{c} n,\rho,\beta,c \\
n',\rho',\beta',c'\end{array}\right)$ as a 2-vector space is
$\{(H_1,\Phi_1),\ldots,(H_N,\Phi_N)\}$ with
$(H_{\lambda},\Phi_{\lambda})$ any intertwiner whose matrix of
ranks $R^{(\lambda)}$ is given by
\begin{equation} \label{matriu_rangs_base}
r^{(\lambda)}_{i'i}=\left\{\begin{array}{ll} 1,&\mbox{if
$(i',i)\in X_{\lambda}$} \\ 0,&\mbox{otherwise}\end{array}\right.
\end{equation}
On the other hand, a morphism
$(f_1,\ldots,f_N):(k^{d_1},\ldots,k^{d_N})\To(k^{\tilde{d}_1},\ldots,k^{\tilde{d}_N})$
gets mapped to the unique 2-intertwiner
$\tau:H\Rightarrow\tilde{H}$ whose components on the
basis~\footnote{The reader may think of $\Vv=\ev_k^n$, in which
case this basis is given by the objects
$(0,\ldots,\stackrel{i)}{k},\ldots,0)$ for all $i=1,\ldots,n$.}
$\{V_1,\ldots,V_n\}$ of $\Vv$ are given by the matrices
$M_1,\ldots,M_N$ of the linear maps $f_1,\ldots,f_N$ in canonical
bases. More precisely, if $(i',i)\in X_{\lambda}$ the restriction
of the map
$\tau_{V_i}:\oplus_{i'}r_{i'i}V'_{i'}\To\oplus_{i'}\tilde{r}_{i'i}V'_{i'}$
to its $V'_{i'}$-`isotypic' component is that defined by the
$(\tilde{d}_{\lambda}\times d_{\lambda})$-matrix $M_{\lambda}$.
Let us emphasize that different matrices $M_{\lambda}$ are used to
define the same morphism $\tau_{V_i}$, and that there is no
obvious general relation between the number of these basic
components $\tau_{V_i}$, which is equal to the rank $n$ of $\Vv$,
and the number of matrices we use to compute them, which is equal
to the number $N$ of intertwining orbits. The equivalence so
defined is clearly non-canonical. It depends, among other things,
on the linear order chosen in the set of intertwining orbits. }
\end{rem}

\subsection{Intertwining numbers}
Let us suppose that $k$ is of characteristic zero and that $\G$ is
essentially finite. In particular, all hom-categories
$\Hh\left(\begin{array}{c} n,\rho,\beta,c \\
n',\rho',\beta',c'\end{array}\right)$ are 2-vector spaces. Then it
follows from the discussion in
\S~\ref{representacions_projectives} that the ranks of these
2-vector spaces or {\it intertwining numbers} can be explicitly
computed by the following procedure:

\begin{itemize}
\item Find the intertwining $\pi_0(\G)$-orbits $X_1,\ldots,X_N$ of $X(n',n)$.
\item For each $\lambda=1,\ldots,N$ choose any point $(i'_{\lambda},i_{\lambda})\in X_{\lambda}$, determine the corresponding stabilizer $G_{i'_{\lambda},i_{\lambda}}\subseteq\pi_0(\G)$, and compute the above normalized 2-cocycle $\hat{z}_{\lambda}:G_{\lambda}\times G_{\lambda}\To k^*$.
\item For each $\lambda=1,\ldots,N$ compute the number $r(G_{\lambda},\hat{z}_{\lambda})$ of $\hat{z}_{\lambda}$-regular conjugacy classes of $G_{\lambda}$.
\end{itemize}
Then the intertwining number is given by
$$
{\rm rank}\ \Hh\left(\begin{array}{c} n,\rho,\beta,c \\ n',\rho',\beta',c'\end{array}\right)=\sum_{\lambda=1}^Nr(G_{\lambda},\hat{z}_{\lambda}).
$$

\vspace{0.2 truecm} \noindent Observe that, proceeding in this
way, we only need to take into account that the morphism is from
$(n,\rho,\beta,c)$ to $(n',\rho',\beta',c')$, and not in the
reverse direction, when computing the 2-cocycles
$\hat{z}_{\lambda}$. However, reversing the direction just
corresponds to replacing $\hat{z}_{\lambda}$ by the inverse
2-cocycle $\hat{z}_{\lambda}^{-1}$. Since the regularity condition
of an element $g\in G_{\lambda}$ is the same either with respect
to $\hat{z}_{\lambda}$ or with respect to
$\hat{z}^{-1}_{\lambda}$, it follows that
$r(G_{\lambda},\hat{z}_{\lambda})=r(G_{\lambda},\hat{z}^{-1}_{\lambda})$.
Hence we have the following analog of the well known symmetry
property for the intertwining numbers between linear
representations of a finite group.

\begin{cor}
Under the above assumptions on $\G$ and $k$ we have
$$
{\rm rank}\ \Hh\left(\begin{array}{c} n,\rho,\beta,c \\ n',\rho',\beta',c'\end{array}\right)={\rm rank}\ \Hh\left(\begin{array}{c} n',\rho',\beta',c' \\ n,\rho,\beta,c\end{array}\right)
$$
for any quadruples $(n,\rho,\beta,c),(n',\rho',\beta',n')$.
\end{cor}

\begin{ex}{\rm
Let us think of the symmetric group $S_n$ as the group of
automorphisms of the finite set $[n]$. Then for any morphism of
groups $\rho:\pi_0(\G)\To S_n$ we have
$$
{\rm rank}\ \Hh\left(\begin{array}{c} 1,1,1,1 \\ n,\rho,1,1\end{array}\right)=\sum_{\lambda=1}^s|\{\mbox{conjugacy classes of $G_s$}\}|,
$$
where $G_1,\ldots,G_s$ denote the stabilizers (determined up to
conjugation) of the $\pi_0(\G)$-orbits $X_1,\ldots,X_s$ of $[n]$.
In particular, if we take as $\rho$ Cayley's morphism
$\rho_C:\pi_0(\G)\To S_{|\pi_0(\G)|}$ we obtain
$$
{\rm rank}\ \Hh\left(\begin{array}{c} 1,1,1,1\\
|\pi_0(\G)|,\rho_C,1,1\end{array}\right)= |\{\mbox{conjugacy
classes of $\pi_0(\G)$}\}|,
$$
in agreement with the fact that we have an equivalence of
$k$-additive categories
$$
\Hh\left(\begin{array}{c} 1,1,1,1\\
|\pi_0(\G)|,\rho_C,1,1\end{array}\right)\simeq\rep_{\ev_k}(\pi_0(\G))
$$
(see Theorem~\ref{hom-categories}). In particular, this is true
when $n=1$ and hence, for the category
$\mathcal{E}nd_{\G}(\Ii)=\Hh\left(\begin{array}{c} 1,1,1,1 \\
1,1,1,1\end{array}\right)$. }
\end{ex}

\begin{ex}\label{exemple2}{\rm
Let $(n,\rho,\beta,c)$ be any quadruple of the kind described in
\S~\ref{classes_iso_repr} and let
$(n_{\Rb},\rho_{\Rb},\beta_{\Rb},c_{\Rb})$ be the quadruple
classifying the regular representation of $\G$ (see
\S~\ref{classificacio}). In particular, we know that $n_{\Rb}=pq$.
Then we have
\begin{equation}\label{rang}
{\rm rank}\ \Hh\left(\begin{array}{c} n_{\Rb},\rho_{\Rb},\beta_{\Rb},c_{\Rb}\\ n,\rho,\beta,c\end{array}\right)=n.
\end{equation}
To see this, note first that the stabilizer of any point
$(i,k,l)\in [n]\times[q]\times[p]\cong X(n,pq)$ is trivial.
Indeed, the point $(k,l)\in[q]\times[p]$ corresponds to the pair
$(\chi_k,g_l)\in\pi_1(\G)^*\times\pi_0(\G)$ (we work with the
linear orders we have fixed in \S~\ref{classificacio} for
$\pi_0(\G)$ and $\pi_1(\G)^*$) . Hence the action of
$g_{l'}\in\pi_0(\G)$ on $(i,k,l)$ is
$$
(i,k,l)\cdot g_{l'}=(i,\chi_k,g_l)\cdot
g_{l'}=(\rho(g^{-1}_{l'})(i),\chi_k,g_lg_{l'})
$$
(cf. Proposition~\ref{classificacio_repr_regular}), and
$$
(\rho(g^{-1}_{l'})(i),\chi_k,g_lg_{l'})\neq(i,\chi_k,g_l)=(i,k,l)
$$
unless $g_{l'}=e$. It follows from Corollary~\ref{corol.lari} that
\begin{equation} \label{homs_R_F}
\Hh\left(\begin{array}{c} n_{\Rb},\rho_{\Rb},\beta_{\Rb},c_{\Rb}\\ n,\rho,\beta,c\end{array}\right)\simeq\ev_k^{N},
\end{equation}
where $N=N(n_{\Rb},\rho_{\Rb},\beta_{\Rb};n,\rho,\beta)$ is the
number of intertwining orbits for the given pair of
representations. It remains to see that $N=n$. In fact, we shall
determine explicitly the intertwining orbits by identifying a
`canonical' representative point in each of them. Let us fix a
character $\chi\in\pi_1(\G)^*$ and a $\pi_0(\G)$-orbit $\Oo$ of
$[n]_{\rho}$. The subset
$$
X_{\Oo,\chi}:=\Oo\times\{\chi\}\times\pi_0(\G)\subset X(n,pq)
$$
is $\pi_0(\G)$-invariant but non-transitive. For example, for any
$i\neq i'$ in $\Oo$ the points $(i,\chi,e)$ and $(i',\chi,e)$ are
not in the same orbit. Actually, the set $\{\ (i,\chi,e),\
i\in\Oo\}$ constitutes a set of representative points for the
various orbits of $X_{\Oo,\chi}$. Indeed, an arbitrary point
$(i,\chi,g)\in X_{\Oo,\chi}$ is in the same orbit as
$(\rho(g)(i),\chi,e)$. Therefore, the decomposition of
$X_{\Oo,\chi}$ into orbits looks like
$$
X_{\Oo,\chi}=\coprod_{i\in\Oo}\ X_{i,\chi},
$$
with $X_{i,\chi}:=(i,\chi,e)\pi_0(\G)$. In particular,
$X_{\Oo,\chi}$ has $|\Oo|$ orbits, all of them of cardinal $p$.
Since this is true for each $\chi\in\pi_1(\G)^*$ it follows that
the decomposition of $X(n,pq)$ into orbits is
$$
X(n,pq)=\coprod_{\chi\in\pi_1(\G)^*}\ \coprod_{\Oo\in{\rm Orb}([n]_{\rho})}\ \coprod_{i\in\Oo}\ X_{i,\chi}=\coprod_{\chi\in\pi_1(\G)^*}\ \coprod_{i\in[n]}X_{i,\chi}.
$$
However, only $n$ of these $qn$ orbits are intertwining. This is
because the $(\chi,e)$-component of $\beta_{\Rb}$ is
$\beta_{\Rb,(\chi,e)}=\chi$ (see
Proposition~\ref{classificacio_repr_regular}). Thus $(i,\chi,e)\in
X_{i,\chi}$ is intertwining if and only if $\beta_i=\chi$.
Consequently the set of intertwining orbits is
\begin{equation}\label{orbites_intertwining}
\Lambda(n_{\Rb},\rho_{\Rb},\beta_{\Rb};n,\rho,\beta)=\{X_{i,\beta_i},i\in[n]\}=\{(i,\beta_i,e)\pi_0(\G),\
i\in[n]\},
\end{equation}
with $\{(i,\beta_i,e),\ i=1,\ldots,n\}$ as a set of `canonical'
representatives. In particular $N=n$ as claimed. }
\end{ex}
Once more, this example is nothing but the analog in our setting of
a similar result concerning the regular representation of a finite
group $G$. Actually, as in the group setting, this
is one of the consequences of the more fundamental fact that for any essentially
finite 2-group $\G$ the forgetful 2-functor
$\mbox{\boldmath$\omega$}:\repg_{\dev_k}(\G)\To\dev_k$ is
represented by the regular representation (se Section~\ref{representabilitat}).

\subsection{Remarks about $k$-linear enrichments on 2-categories}
By a $k$-{\it linear} (resp. $k$-{\it additive}) 2-category we
mean a 2-category $\Cgg$ such that all its hom-categories
$\homs_{\Cgg}(X,Y)$ are $k$-linear (resp. $k$-additive) and all
composition functors
$$
\homs_{\Cgg}(X,Y)\times\homs_{\Cgg}(Y,Z)\To\homs_{\Cgg}(X,Z)
$$
are $k$-bilinear. More particularly, a 2-category will be called a
$\dev_k$-{\it category}~\footnote{It would be more correct to
speak of $2\ev_k$-categories, where $2\ev_k$ stands for the
underlying category of $\dev_k$ equipped with the monoidal
structure defined by the tensor product of $k$-additive categories
(2-vector spaces are indeed stable under this tensor product).}
when it is $k$-additive and all its
hom-categories are 2-vector spaces.

The simplest example of a $\dev_k$-category is $\dev_k$ itself,
which is supposed to be (monoidal) {\it pseudo-closed} in the
sense of \cite{HP02}. Another example is $\repg_{\dev_k}(\G)$ for
$\G$ an essentially finite 2-group and $k$ an algebraically closed
field of characteristic zero or prime to the order of $\pi_0(\G)$.
Indeed, we have shown that all hom-categories in
$\repg_{\dev_k}(\G)$ are 2-vector spaces under these assumptions,
and the corresponding composition functors are $k$-bilinear
because they are so in $\dev_k$. For arbitrary $\G$ and $k$,
$\repg_{\dev_k}(\G)$ is just a $k$-additive 2-category (cf.
\S~\ref{homs_G}).

Given two $k$-linear 2-categories $\Cgg$ and $\Dgg$, a
pseudofunctor $\Fb:\Cgg\To\Dgg$ is called $k$-linear when all
functors
$\Fb_{X,Y}:\homs_{\Cgg}(X,Y)\To\homs_{\Dgg}(\Fb(X),\Fb(Y))$ are
$k$-linear. An example is the forgetful 2-functor
$\mbox{\boldmath$\omega$}:\repg_{\dev_k}(\G)\To\dev_k$ mentioned
before.

Let us finally remark that for any $k$-linear 2-category $\Dgg$
the 2-category ${\bf PsFun}(\Cgg,\Dgg)$ of pseudofunctors from any
other 2-category $\Cgg$ to $\Dgg$ is also $k$-linear, as the
reader may easily check. However, this fails to be true when
$k$-linear is replaced by $k$-additive. For instance, there may
exist no zero object in the hom-categories of ${\bf
PsFun}(\Cgg,\Dgg)$ even when we have a zero object in each
hom-category $\homs_{\Dgg}(X,Y)$. However, we are interested in
cases where $\Dgg$ is $\dev_K$, and ${\bf PsFun}(\Cgg,\dev_k)$ is
always $k$-additive. This is because the objects in $\dev_k$ are
themselves categories with a zero object and binary biproducts,
and these can be used to get a zero object and binary biproducts
in ${\bf PsFun}(\Cgg,\dev_k)$. This is in fact how we have seen
before that $\repg_{\dev_k}(\G)$ is $k$-additive.
The same thing works for the 2-category of pseudofunctors between
$\repg_{\dev_k}(\G)$ and $\dev_k$. In particular, the category
$\mathcal{E}nd(\mbox{\boldmath$\omega$})$ of (pseudonatural)
endomorphisms of $\mbox{\boldmath$\omega$}$ is always
$k$-additive. We shall see in the next section that it is even a
2-vector space under suitable assumptions.

\section{Representability of the forgetful 2-functor $\mbox{\boldmath$\omega$}$}
\label{representabilitat}

In order to prove the representability of
$\mbox{\boldmath$\omega$}$ we shall make use of the appropriate
enriched version of the bicategorical Yoneda Lemma. Hence this
section starts by recalling this basic result as well as the
associated notion of ``universal object'' for $\catg$-valued
($\dev_k$-valued in the enriched case) pseudofunctors. These are analogs
of the universal elements of a $\sets$-valued (resp.
$\ev_k$-valued) functor. Next it is shown that for essentially
finite 2-groups $\G$ and algebraically closed fields $k$ as above
there indeed exists a universal object for the forgetful 2-functor
$\mbox{\boldmath$\omega$}:\repg_{\dev_k}(\G)\To\dev_k$ leading to
a representation of it by the regular representation of $\G$. The
section closes with a description of this representation and the
induced equivalence between the category
$\mathcal{E}nd(\mbox{\boldmath$\omega$})$ of pseudonatural
endomorphisms of $\mbox{\boldmath$\omega$}$ and $\ev_k^{\Gg}$. As
mentioned in the introduction, this equivalence constitutes a
first step toward a Tannaka-Krein type reconstruction of an
essentially finite 2-group from its 2-category of representations
in 2-vector spaces (and the associated forgetful 2-functor).

\subsection{Bicategorical Yoneda Lemma}
To my knowledge, this result first appears in its nonenriched version in
\cite{rS80}. It establishes the (natural) 
equivalence of two categories. More specifically, suppose given a
bicategory $\Cgg$, an object $X$ of $\Cgg$, and a pseudofunctor
$\Ff:\Cgg\To\catg$ with values in the 2-category $\catg$ of
(small) categories, functors and natural transformations. Let
us denote by $\homs_{\Cgg}(X,-):\Cgg\To\catg$ the (covariant)
hom-pseudofunctor associated to $X$, and let
$\mathcal{P}sNat(\homs_{\Cgg}(X,-),\Ff)$ be the category of all
pseudonatural transformations $\homs_{\Cgg}(X,-)\Rightarrow\Ff$
and modifications between these. Then the bicategorical Yoneda
lemma says that there exists an equivalence (not an isomorphism)
of categories
$$
{\sf
Yon}:\mathcal{P}s\mathcal{N}at(\homs_{\Cgg}(X,-),\Ff)\stackrel{\simeq}{\To}\Ff(X)
$$
which is given on objects $\xi:\homs_{\Cgg}(X,-)\Rightarrow\Ff$ by
$$
{\sf Yon}(\xi):=\xi_X({\id}_X),
$$
and that this equivalence is natural in $X$ and $\Ff$ in some
suitable sense. Unlike usually for equivalences of categories
which are not isomorphisms, {\sf Yon} has a canonical
pseudoinverse. In fact, although we shall make no use of it, it
can be shown that $\yon$ extends canonically to an adjoint
equivalence $({\sf Yon},{\sf Yon}^*,\eta,\epsilon)$ whose unit
$\eta$ is an identity when $\Cgg$ is a 2-category, and whose
counit $\epsilon$ is an identity when $\Ff$ is a (strict)
2-functor. The canonical pseudoinverse
$$
{\sf
Yon}^*:\Ff(X)\To\mathcal{P}s\mathcal{N}at(\homs_{\Cgg}(X,-),\Ff)
$$
maps $A\in{\rm Obj}\Ff(X)$ to the pseudonatural transformation
${\sf Yon}^*(A):\homs_{\Cgg}(X,-)\Rightarrow\Ff$ whose 1-cell
components ${\sf Yon}^*(A)_Y:\homs_{\Cgg}(X,Y)\To\Ff(Y)$ are given
on objects $f:X\To Y$ and morpisms $\tau:F\Rightarrow f'$ by
\begin{align}
{\sf Yon}^*(A)_Y(f)&:=\Ff(f)(A) \label{yon*1} \\ {\sf
Yon}^*(A)_Y(\tau)&:=\Ff(\tau)_A. \label{yon*2}
\end{align}
What we really need is the $k$-linear version of this lemma. In
this version $\Cgg$ is a $k$-linear 2-category,
$\Ff:\Cgg\To\dev_k$ a $k$-linear pseudofunctor and
\begin{equation} \label{equivalencia_Yoneda}
{\sf
Yon}:\mathcal{P}s\mathcal{N}at(\homs_{\Cgg}(X,-),\Ff)\stackrel{\simeq}{\To}\Ff(X)
\end{equation}
a $k$-linear equivalence and hence, an equivalence of 2-vector
spaces. In particular, the $k$-linear category
$\mathcal{P}s\mathcal{N}at(\homs_{\Cgg}(X,-),\Ff)$ is actually a
2-vector space. The pseudoinverse $\yon^*$ is defined in exactly
the same way as before.

\subsection{Universal objects of a pseudofunctor}
Given a representable functor $F:\Cc\To\sets$ with representing
object $X\in{\rm Obj}\ \Cc$, it is well known that the Yoneda
bijection
$$
Yon:{\rm Nat}({\rm Hom}_{\Cc}(X,-),F)\stackrel{\cong}{\To} F(X)
$$
restricts to a bijection between representations of $F$ by $X$ (isomorphisms
${\rm Hom}_{\Cc}(X,-)\cong F$) and the so called
{\it universal elements} in $F(X)$. These are elements $x\in
F(X)$ such that for any object $Y$ of $\Cc$ the map ${\rm
Hom}_{\Cc}(X,Y)\To F(Y)$ defined by $f\mapsto F(f)(x)$ is a
bijection. The result follows directly from the definition of the
bijection to $Yon$. The analogous result holds
for $\ev_k$-valued functors and the corresponding Yoneda
isomorphisms of vector spaces.

Similarly, given any $\catg$-valued or $\dev_k$-valued
pseudofunctor $\Ff$ on a 2-category $\Cgg$ and any object $X$ of
$\Cgg$, by a {\it universal object} of $\Ff(X)$ we mean an object
$x\in{\rm Obj}\Ff(X)$ such that the pseudonatural transformation
$\yon^*(x)$ is a pseudonatural equivalence. Now it is a general
fact that a pseudonatural transformation is an equivalence if and
only if all its 1-morphism components are equivalences. Hence
$x\in{\rm Obj}\Ff(X)$ is universal if and only if the functors
$$
\yon^*(x)_Y:\homs_{\Cgg}(X,Y)\To\Ff(Y)
$$
are equivalences of categories for any $Y\in{\rm Obj}\Cgg$.

By the very definition of universal objects, it follows that the
Yoneda equivalence (\ref{equivalencia_Yoneda}) restricts to a
($k$-linear) equivalence of categories
$$
\yon:\pseq(\homs_{\Cgg}(X,-),\Ff)\stackrel{\simeq}{\To}\Ff(X)_u
$$
between the full subcategory $\pseq(\homs_{\Cgg}(X,-),\Ff)$ with
objects all pseudonatural equivalences (i.e. 
representations of $\Ff$ by $X$) and the full subcategory
$\Ff(X)_u$ with objects the universal ones. In particular,
the pseudofunctor $\Ff$ is representable by the object $X$ of
$\Cgg$ or equivalently,
$\mathcal{P}s\mathcal{E}q(\homs_{\Cgg}(X,-),\Ff)$ is nonempty if
and only if there exists such a universal object $x\in{\rm
Obj}\Ff(X)$.

\subsection{Universal functor $\eta_U:\Gg\To\ev_k$}
\label{seccio_functor_universal} We are interested in the case
where $\Cgg$ is the 2-category $\repg_{\dev_k}(\G)$ and $\Ff$ the
forgetful 2-functor
$\mbox{\boldmath$\omega$}:\repg_{\dev_k}(\G)\To\dev_k$. According
to the previous discussion, in order to prove that
$\mbox{\boldmath$\omega$}$ is represented by the regular
representation $\Rb$ it is enough to see that there exists a
universal object in the 2-vector space
$\mbox{\boldmath$\omega$}(\Rb)=\ev_k^{\Gg}$. More precisely, we
have to see that there exists a functor
$$
\eta_U:\Gg\To\ev_k
$$
satisfying the next two conditions:
\begin{itemize}
\item[(i)] For any representation $\Fb=(\Vv,\F)$ and any $V\in{\rm Obj}\ \Vv$ there exists an intertwiner $(H,\Phi):\Rb\To\Fb$ such that $H(\eta_U)\cong V$ (i.e. $\yon^*(\eta_U)_{\Fb}$ is essentially
surjective; cf. (\ref{yon*1})).
\item[(ii)] For any representation $\Fb=(\Vv,\F)$, any intertwiners $(H,\Phi),(\tilde{H},\tilde{\Phi}):\Rb\To\Fb$ and any morphism $\phi:H(\eta_U)\To\tilde{H}(\eta_U)$ in $\Vv$ there exist a unique 2-intertwiner $\tau:(H,\Phi)\Rightarrow(\tilde{H},\tilde{\Phi})$ such that
$\phi=\tau_{\eta_U}$  (i.e. $\yon^*(\eta_U)_{\Fb}$ is fully
faithful; cf. (\ref{yon*2})).
\end{itemize}
We claim that such a universal functor exists and is given by the
direct sum of a few of the basic functors $\{\eta_{\chi,g},\
(\chi,g)\in\pi_1(\G)^*\times\pi_0(\G)\}$ of Example~\ref{vectG}. More
explicitly, we have the following. 

\begin{thm}\label{functor_universal}
Let $\G$ be essentially finite and $k$ algebraically closed and of
characteristic zero. Then the pair $(\Rb,\eta_U)$, with $\Rb$ the
regular representation of $\G$ and
$$
\eta_U:=\bigoplus_{\chi\in\pi_1(\G)^*}\eta_{\chi,e}:\Gg\To\ev_k,
$$
is a universal object for the forgetful 2-functor
$\mbox{\boldmath$\omega$}:\repg_{\dev_k}(\G)\To\dev_k$. In
particular, $\mbox{\boldmath$\omega$}$ is representable with $\Rb$
as representing object.
\end{thm}
\begin{proof}
Suppose we are given $\Fb$ and $V$ as in (i). Let
$(n,\rho,\beta,c)$ be the quadruple classifying $\Fb$, and let
$\{V_1,\ldots,V_n\}$ be a basis of absolutely simple objects of
$\Vv$. We know that, up to a 2-isomorphism, any intertwiner
$(H,\Phi):\Rb\To\Fb$ is completely determined by the matrix of
ranks of $H:\ev_k^{\Gg}\To\Vv$ (see (\ref{homs_R_F}) and
Remark~\ref{remarca_matriu_rangs}). Let $R=(r_{j,(\chi,g)})$ be
this matrix. Thus we have
$$
H(\eta_{\chi,g})\cong\bigoplus_{j=1}^nr_{j,(\chi,g)}V_j
$$
for any $(\chi,g)\in\pi_1(\G)^*\times\pi_0(\G)$. It follows from
the invariance properties of $R$ and the fact that it is supported
on the intertwining orbits that this matrix is actually completely
given by $n$ integers $d_1,\ldots,d_n\geq 0$ giving the values of
the nonzero ``intertwining entries''. To be precise, we shall
assume that $d_i$ gives the value of the entries of $R$
corresponding to the intertwining orbit $X_{i,\beta_i}$ (see
(\ref{orbites_intertwining}) for notation). Then let us take as
$(H,\Phi)$ any intertwiner whose isomorphism class is determined
in this way by the unique integers $d_1,\ldots,d_n\geq 0$ such
that
$$
V\cong\bigoplus_{i=1}^nd_iV_i.
$$
Thus $H$ is a $k$-linear functor whose matrix of ranks is given by
$$
r_{j,(\chi,g)}=\left\{\begin{array}{ll} d_i,& \mbox{if
$(j,\chi,g)\in X_{i,\beta_i}$ for some $i\in\{1,\ldots,n\}$} \\
0,& \mbox{otherwise} \end{array}\right. .
$$
Now, the action of $\pi_0(\G)$ on $X(n,pq)\cong[n]\times\pi_1(\G)^*\times\pi_0(\G)$ is given by
$$
(i,\chi,g)\cdot
\tilde{g}=(\rho(\tilde{g}^{-1})(i),\chi,g\tilde{g}^{-1})
$$
(see \S~\ref{classificacio}). Since
$X_{i,\beta_i}=(i,\beta_i,e)\pi_0(\G)$ it follows that
\begin{align*}
r_{j,(\chi,g)}\neq 0\ &\Leftrightarrow\ \exists\
i\in\{1,\ldots,n\},\ \exists\ \tilde{g}\in\pi_0(\G):
(j,\chi,g)=(\rho(\tilde{g}^{-1})(i),\beta_i,\tilde{g}^{-1})
\\ &\Leftrightarrow\ \exists\ i\in\{1,\ldots,n\}:\ j=\rho(g^{-1})(i),\ \chi=\beta_i,
\end{align*}
in which case we have $r_{j,(\chi,g)}=d_i=d_{\rho(g)(j)}$. If we
define
$$
J(\chi,g):=\{j\in\{1,\ldots,n\}\ |\ \exists\ i\in\{1,\ldots,n\}:\
\rho(g^{-1})(i)=j,\ \beta_i=\chi\}
$$
it follows that
$$
H(\eta_{\chi,g})\cong\bigoplus_{j\in J(\chi,g)}d_{\rho(g)(j)}V_j.
$$
In particular we have
\begin{equation} \label{eta_chi_e}
H(\eta_{\chi,e})\cong\bigoplus_{j\in
J(\chi,e)}d_jV_j=\bigoplus_{\beta_j=\chi}d_jV_j,
\end{equation}
where the condition $\chi=\beta_j$ in the last direct sum means
that it is taken over all $j\in\{1,\ldots,n\}$ such that
$\beta_j=\chi$. Using now that $H$ is $k$-linear we obtain that
\begin{equation} \label{H_eta_U}
H(\eta_U)\cong\bigoplus_{\chi\in\pi_1(\G)^*}H(\eta_{\chi,e})\cong\bigoplus_{\chi\in\pi_1(\G)^*}\left(\bigoplus_{\beta_j=\chi}d_jV_j\right)\cong\bigoplus_{j=1}^nd_jV_j\cong V.
\end{equation}
This proves (i).

To prove (ii) let us first remark that for any intertwiners
$(H,\Phi)$, $(\tilde{H},\tilde{\Phi})$ from $\Rb$ to any other
representation $\Fb$ we have a bijection
$$
A:{\rm
2Hom}_{\G}((H,\Phi),(\tilde{H},\tilde{\Phi}))\stackrel{\cong}{\longrightarrow}\prod_{i=1}^n{\rm
Hom}_k(k^{d_i},k^{\tilde{d}_i})
$$
between the set of 2-intertwiners
$\tau:(H,\Phi)\Rightarrow(\tilde{H},\tilde{\Phi})$, on the one
hand, and the set of $n$-tuples of linear maps $(f_1,\ldots,f_n)$
with $f_i:k^{d_i}\To k^{\tilde{d}_i}$, on the other. It basically give the
morphisms of vector bundles over the various 
intertwining orbits $X_{i,\beta_i}$, and its existence follows from
(\ref{homs_R_F}). Moreover, it follows from Remark~\ref{remarca_matriu_rangs}
that this bijection maps a 2-intertwiner $\tau$ to the $n$ linear maps
$A(\tau)_1,\ldots,A(\tau)_n$ obtained in the following way. From
(\ref{eta_chi_e}) we know that
$$
H(\eta_{\beta_i,e})\cong\bigoplus_{j=1}^nr_{j,(\beta_i,e)}V_j=\bigoplus_{\beta_j=\beta_i}d_jV_j,
$$
where the last direct sum is taken over all $j\in\{1,\ldots,n\}$
such that $\beta_j=\beta_i$, and similarly for
$\tilde{H}(\eta_{\beta_i,e})$. Hence we have
$$
\tau_{\eta_{\beta_i,e}}:\bigoplus_{\beta_j=\beta_i}d_jV_j\To\bigoplus_{\beta_j=\beta_i}\tilde{d}_jV_j,
$$
Then $A(\tau)_i:k^{d_i}\To k^{\tilde{d}_i}$ is the linear map
whose matrix in canonical bases is equal to the matrix giving the
restriction of $\tau_{\eta_{\beta_i,e}}$ to the $V_i$-`isotypic'
component.

Suppose now we are given a morphism
$\phi:H(\eta_U)\To\tilde{H}(\eta_U)$ in $\Vv$. Because of
(\ref{H_eta_U}) and the absolute simplicity of the objects $V_i$,
we see that giving such a morphism amounts to giving $n$
arbitrary linear maps $f_i:k^{d_i}\To k^{\tilde{d}_i}$. Therefore
we also have a bijection
$$
B:\prod_{i=1}^n{\rm
Hom}_k(k^{d_i},k^{\tilde{d}_i})\stackrel{\cong}{\longrightarrow}{\rm
Hom}_{\Vv}(H(\eta_U),\tilde{H}(\eta_U)),
$$
mapping an $n$-tuple $(f_1,\ldots,f_n)$ to the
unique morphism $\phi:H(\eta_U)\To\tilde{H}(\eta_U)$ whose
restriction to the $V_i$-`isotypic' component is given by the
matrix of $f_i$ in canonical bases. It then follows that the
composite bijection
$$
B\circ A:{\rm
2Hom}_{\G}((H,\Phi),(\tilde{H},\tilde{\Phi}))\longrightarrow{\rm
Hom}_{\Vv}(H(\eta_U),\tilde{H}(\eta_U))
$$
maps any $\tau$ to its component $\tau_{\eta_U}$, and this proves
(ii).
\end{proof}

\begin{cor}
Under the same assumptions on $\G$ and $k$ as before we have an
equivalence of $k$-additive categories
$$
\mathcal{E}nd(\mbox{\boldmath$\omega$})\simeq\ev_k^{\Gg}.
$$
In particular, $\mathcal{E}nd(\mbox{\boldmath$\omega$})$ is a
2-vector space of rank $pq$.
\end{cor}
\begin{proof}
Any equivalence $f:X\stackrel{\simeq}{\To} Y$ in a ($k$-linear)
2-category $\Cgg$ induces equivalences of ($k$-linear) categories
$$
\mathcal{E}nd_{\Cgg}(X)\simeq\homs_{\Cgg}(X,Y)\simeq\mathcal{E}nd_{\Cgg}(Y).
$$
In our case $\Cgg$ is the 2-category ${\bf
PsFun}(\repg_{\dev_k}(\G),\dev_k)$, $X$ is the hom-pseudofunctor
$\homs_{\G}(\Rb,-)$, $Y$ is $\mbox{\boldmath$\omega$}$ and $f$ is
any representation of $\mbox{\boldmath$\omega$}$ by $\Rb$. Hence
we have
$$
\mathcal{E}nd(\mbox{\boldmath$\omega$})\simeq\mathcal{P}s\mathcal{N}at(\homs_{\G}(\Rb,-),\mbox{\boldmath$\omega$})
\simeq\ev_k^{\Gg}
$$
because of Yoneda.
\end{proof}

\begin{rem}{\rm
It is worth comparing this with the situation we have for finite
groups. Thus for any finite group $G$ there also exists a
``universal function'' $f_U:G\To k$, i.e. a function such that for
any representation $V$ and any $v\in V$ there exists a unique
morphism of representations $h:L(G)\To V$ with $h(f_U)=v$. Such a
universal function is the function $\delta_e$ equal to zero
everywhere except on the unit element $e\in G$ where it is equal
to 1. Hence the analog in our categorified setting of the basic
function $\delta_e$ is none of the basic functors $\eta_{\chi,e}$,
for some particular $\chi\in\pi_1(\G)^*$, but the direct sum of
all of them.

}
\end{rem}

\subsection{Basis of
  $\mathcal{E}nd(\mbox{\boldmath$\omega$})$}

From Theorem~\ref{functor_universal} and the above
description of ${\sf Yon}^*$ (see (\ref{yon*1})-(\ref{yon*2})) it follows 
that a specific representation of
$\mbox{\boldmath$\omega$}:\repg_{\dev_k}(\G)\To\dev_k$ is the 
2-natural equivalence
$$
\Theta\equiv\yon^*(\eta_U):\homs_{\G}(\Rb,-)\Rightarrow\mbox{\boldmath$\omega$}
$$
whose 1-cell components are the $k$-linear functors
$\Theta_{\Fb}:\homs_{\G}(\Rb,\Fb)\To\Vv$ given on objects $(H,\Phi)$ and
morphisms $\tau:H\Rightarrow H'$ by
\begin{align*}
\Theta_{\Fb}(H,\Phi)&:=H(\oplus_{\chi}\eta_{\chi,e})
\\ \Theta_{\Fb}(\tau)&:=\tau_{\oplus_{\chi}\eta_{\chi,e}}.
\end{align*}
This induces a ($k$-linear) equivalence
$$
\mathcal{E}nd(\mbox{\boldmath$\omega$})\To\mathcal{P}s\mathcal{N}at(\homs_{\G}(\Rb,-)\mbox{\boldmath$\omega$})
$$
defined by $u\mapsto u\cdot\Theta$. Composing it with the Yoneda
equivalence gives the desired $k$-linear equivalence
$E:\mathcal{E}nd(\mbox{\boldmath$\omega$})\stackrel{\simeq}{\longrightarrow}\ev_k^{\Gg}$. This
turns out to be the equivalence mapping the pseudonatural transformation
$u:\mbox{\boldmath$\omega$}\Rightarrow\mbox{\boldmath$\omega$}$ to
the functor
$$
u_{\Rb}(\oplus_{\chi}\eta_{\chi,e}):\Gg\To\ev_k,
$$
and a modification $\ngg:u\Rrightarrow u'$ to the natural
transformation
$$
\ngg_{\Rb,\oplus_{\chi}\eta_{\chi,e}}:u_{\Rb}(\oplus_{\chi}\eta_{\chi,e})\Rightarrow
u'_{\Rb}(\oplus_{\chi}\eta_{\chi,e}).
$$
A pseudoinverse $E^*:\ev_k^{\Gg}\To\mathcal{E}nd(\mbox{\boldmath$\omega$})$ is
given by
$$
E^*(\eta)=\yon^*(\eta)\cdot\Theta^*,\qquad \eta\in{\rm Obj}\ev_k^{\Gg}
$$
for some pseudoinverse $\Theta^*$ of the above 2-natural equivalence
$\Theta$. The 1-cell components of such a $\Theta^*$ are described in
Remark~\ref{remarca_matriu_rangs} for $\Vv$ of the form
$\ev_k^n$. For an arbitrary 2-vector space $\Vv$ we just need to identify the
standard basis of $\ev_k^n$ with a basis of $\Vv$. In particular, for any
representation $\Fgg=(\Vv,\F)$ the 
$k$-linear equivalence
$$
\Theta^*_{\Fgg}:\Vv\to\homs_{\G}(\Rb,\Fgg)
$$
maps a
basis $\{V_1,\ldots,V_n\}$ of $\Vv$ to the basic intertwiners whose matrices of
ranks are given by (\ref{matriu_rangs_base}). Hence we have the following.
\begin{prop}
A basis of $\mathcal{E}nd(\mbox{\boldmath$\omega$})$ as a 2-vector space is
given by a family of endomorphisms 
$$
\{\zeta_{\chi,g}\equiv
E^*(\eta_{\chi,g}):\mbox{\boldmath$\omega$}\Rightarrow\mbox{\boldmath$\omega$},\
(\chi,g)\in\pi_1(\G)^*\times\pi_0(\G)\}
$$
whose 1-cell components are given by
$$
\zeta_{\chi,g;\Fgg}(V_i)=\left\{\begin{array}{ll} V_{\rho(g^{-1})(i)},&\mbox{if
        $\chi=\beta_i$}, \\ \\ 0,&\mbox{otherwise}\end{array}\right.
$$ 
if $\Fgg\simeq\Fgg(n,\rho,\beta,c)$. In particular, 
$\zeta_{\chi,g}$ is totally supported on representations whose
$\beta=(\beta_1,\ldots,\beta_n)$ includes the character $\chi$.
\end{prop}

\section{Final remark}

A well known important property of the regular representation
$L(G)$ of a finite group $G$ is that it is equivalent (for
algebraically closed fields $k$) to the direct sum of all
nonequivalent irreducible representations, each one with a
multiplicity exactly equal to its own dimension. Because of the
similarities we have found until now one might be tempted to think
that the same is true for essentially finite 2-groups. However,
on the one hand, in our setting there may exist 
non-irreducible but indecomposable representations. This fact has
been pointed out in \cite{BBFW09} in the even more general
framework of representations of 2-groups in Yetter's measurable
categories, of which our representation theory is a special case. Hence not
every representation will necessarily decompose as a direct sum of
irreducible ones. On the other hand, as pointed out before, there is no Schur's
Lemma in our representation theory, at least in
its usual form, and such lemma seems to be crucial to prove the above mentioned
result for finite groups. Indeed, if $G$ is a finite group and $k$
an algebraically closed field Schur's lemma implies that for any
irreducible representation $V_i$ of $G$ the dimension of ${\rm
Hom}_G(L(G),V_i)$ is precisely equal to the multiplicity of $V_i$ in $L(G)$.
The result mentioned above follows then because $\omega$ is
represented by $L(G)$ so that we have $V_i\cong{\rm
Hom}_G(L(G),V_i)$. When we move to our setting, we  
still have an equivalence of 2-vector spaces
$\Vv\simeq\homs_{\G}(\ev_k^{\Gg},\Vv)$ for any linear
representation $(\Vv,\F)$ of $\G$. However, it is not at all clear
whether the rank of $\homs_{\G}(\ev_k^{\Gg},\Vv)$, for $(\Vv,\F)$
irreducible, coincides with the `number of copies' of it in the
regular representation.

\bibliographystyle{plain}
\bibliography{repr_2grups_arxiv}

\end{document}